%% file: paper2.tex
\input{header}

\begin{document}
\parindent 0pt

\title[Subsets of $\F$]{Some subsets of set $\F$ in the diluted Hofstadter problem}

\author{Jonathan H.B.\ Deane}
\address{Department of Mathematics, University of Surrey, Guildford, GU2 7XH, UK}
\email{j.deane@surrey.ac.uk}
\author{Guido Gentile}
\address{Dipartimento di Matematica, Universit\`a Roma Tre, Roma, I-00146, Italy}
\email{guido.gentile@uniroma3.it }
\date{}

\input{abstract2.tex}

\maketitle

\section{Introduction}

A previous work~\cite{dilh} considered the problem of the existence of integer sequences $q$ defined recursively by
\begin{equation}
\label{rec}
q(n) = q(n - q(n-1)) + f(n)\qquad\mbox{with}\qquad q(1) = 1,
\end{equation}
for integer $n>1$, and where integer sequence $f$, with $f(1) = 0$, is given.
Motivation for studying this problem comes from the question of the
existence of the Hofstadter $q$-sequence~\cite{geb}, which is defined
by~\eqref{rec} with $f(n) = q(n-q(n-2))$ and with $q(2)= 1$.
The background to the problem is given in the first section
of~\cite{dilh} and we do not repeat this material here.

All sequences in the current paper start from index 1 and we consider $q(i)$ for
$i\leq 0$ to be undefined: hence, in the light of the `nested' term
$q(n-q(n-1))$ in~\eqref{rec}, the claim that sequence $q$ exists is
equivalent to
\begin{equation}
\label{qb}
1\leq q(n)\leq n\qquad \mbox{for all $n\in\mathbb N$.}
\end{equation}
We define $\Q$ as the set of all sequences that obey~\eqref{qb}.
If there is an index $d$ at which $q(d) < 1$ or $q(d) > d$, then sequence $q$ is not
in $\Q$ and $q$ said to die at $d$.

Among other things, it was proved in~\cite{dilh} that a sufficient but
not necessary condition for the existence of $q$ is that $f$ be slow, i.e.,
$f(i+1)-f(i)\in\{0, 1\}$ for all $i\in\mathbb N$.

In the rest of the paper, we investigate $\F$, the set of sequences $f$ for which the
corresponding sequence $q$ exists. After explaining the notation, we show
that there exists a bijection between $\Q$ and $\F$, which turns out to be
useful in several of the proofs. We find sharp bounds that any $f\in\F$
must obey, whereupon it becomes clear that there are sequences obeying these
bounds which are not members of $\F$. These `gaps' are a source of difficulty in the study
of~\eqref{rec}.

We then turn our attention to the study of $\F$, proving both negative and
positive results. We construct a variety of subsets of $\F$ and those found
in Theorems~\ref{fastest} and~\ref{thm:2} are particularly notable.

Following these results, we give an algorithm, the Bounding Algorithm, which relates,
approximately, the bounds
obeyed by a set of sequences $f$ to those of the corresponding sequences $q$. This
is useful as an investigative tool, but also leads to a (computer-assisted) proof 
when used in conjunction with Theorem~\ref{thm:2}.

We finish with some conclusions.

\subsection{Notation}

The first term of all sequences considered has index 1, and hence all
sequences considered here are semi-infinite. With a slight abuse of terminology,
we describe such sequences as `infinite' throughout.

As in~\cite{dilh}, we use the notation $\left(a(n)\right)_{n\in\mathbb N} := a(1), a(2), \ldots$ for an
infinite sequence. We sometimes use this notation in formal contexts, such as
in definitions and lemmas, but where it is clear what is
intended, we just write $a$ to mean the whole sequence. Naturally, when the
$n$-th term is meant, we write $a(n)$.
Where we need to refer explicitly to the first few terms of a sequence, $a(n)$ say, we list the values
starting from $a(1)$ --- for example $a = (1, 1, 3, 2)$ means $a(1) = a(2) = 1$, $a(3) = 3$, $a(4) = 2$.

We write the set of consecutive integers $i, i+1, \ldots j$, with $j\geq
i$, as $\{i:j\}$. Hence $\{i:i\}$ is identical to $\{i\}$.

It is sometimes convenient to write $q'(n) := q(n-q(n-1))$ for all
$n\in\mathbb N$, where we define $q'(1) = 1$. Also, given a sequence $q$
obeying~\eqref{qb}, we write $q' = Q'(q)$.

Given two sets of integers, $G$ and $H$, say, their \textit{difference set}
is defined as $G - H := \{g-h | g\in G, h\in H\}$ and their \textit{sumset}
is defined as $G + H := \{g+h | g\in G, h\in H\}$,
with the cardinality of both sets being counted without multiplicity.

In several places we use Stirling's approximation~\cite{a&s} to estimate $n! = \Gamma(n+1)$ by
\begin{equation}
\label{stirling}
\Gamma(n+1) = \sqrt{2\pi n} \left(\frac{n}{\mathrm{e}}\right)^n\left[1 + O\left(\frac{1}{n}\right)\right].
\end{equation}
We use $\Gamma(n+1)$ instead of $n!$ here since it is defined for $n\in\mathbb R$. In practice, since we typically 
need this to estimate growth rates of functions, with a slight abuse of
notation, we use the approximation for $n!$ 
regardless of whether $n$ is an integer or a positive real number.  \\

\section{Sets $\Q$ and $\F$}

We first define two sets, $\Q$ and $\F$, of infinite integer sequences. Let
$\Q$ be the set of all infinite sequences 
$\left(q(n)\right)_{n\in\mathbb N}$ obeying~\eqref{qb}.

Set $\F$ is defined in relation to $\Q$ by solving~\eqref{rec} for $f(n)$: that is, for
any $q\in\Q$, $f$ is in $\F$ if $f(1) = 0$ and
\begin{equation}
\label{f_of_q}
f(n) = q(n) - q(n-q(n-1))
\quad\mbox{for $n>1$.}
\end{equation}
Let $F(q)$ denote the sequence $f$ computed via~\eqref{f_of_q} for given $q\in\Q$.
We also define $Q(f)$ as taking a sequence $f\in\F$ and computing the corresponding
sequence, $q$, using~\eqref{rec}.

We define $\Q_n$ as the set of the first $n$ terms of the sequences in $\Q$,
and $\F_n$ is defined analogously. By the definition of $\Q$, $\Q_n$
contains $n!$ sequences.

\subsection{The bijection between $\Q$ and $\F$}
Let $n>0$ and define the two sequences $A_n = (a(1), \ldots, a(n))$ and $B_n = (b(1), \ldots, b(n))$.
Then $A_n\neq B_n$ if there exists an integer $k$ obeying
$1\leq k\leq n$, such that $a(k)\neq b(k)$; if no such $k$ exists, then $A_n = B_n$.

This definition is used in the following lemma.
\begin{lemma}
\label{injective}
Fix $n>1$ and let $q = (q(1), \ldots, q(n))$, $\hat{q} = (\hat{q}(1), \ldots, \hat{q}(n))$
with $q, \hat{q}\in\Q_n$ and $q\neq \hat{q}$. Let $f = F(q)$ and $\hat{f} = F(\hat{q})$. Then $f\neq \hat{f}$.
\end{lemma}
\begin{proof}
Let $k$ be the smallest integer such that $q(k)\neq \hat{q}(k)$: since $q\neq \hat{q}$, such an integer must exist.
We use~\eqref{f_of_q} to compute $f := F(q) = \left(f(1), \ldots, f(k-1)\right)$.
At the $i$-th step, we have $f(i) = q(i) - q(i-q(i-1))$. Since
$q\in\Q_n$, each $q(i)$ obeys $1\leq q(i)\leq i$, hence, $1\leq i -
q(i-1)\leq i-1$ and so $1\leq q(i-q(i-1)) \leq i-1$. Therefore, $f(i)-q(i)$ for
each $i$ depends only on the value of $q(j)$ with $1\leq j\leq i-1$.

Making the same argument for $\hat{q}$ leads us to conclude that
$\hat{f}(i)-\hat{q}(i)$ depends only on the values of $\hat{q}(j)$ for $1\leq j\leq i-1$.

Now, consider the case $i=k$. We have that $f(k)-q(k)$ depends only on $q(j)$
with $1\leq j\leq k-1$, and so $f(k) - q(k) = \hat{f}(k) - \hat{q}(k)$. Now, since
$q(k) \neq \hat{q}(k)$, we conclude that $f(k)\neq \hat{f}(k)$ and hence, if $q\neq \hat{q}$, $F(q)\neq F(\hat{q})$.
\end{proof}
\begin{remark}
Therefore,
$F \! : \Q_n \to\F_n$ is injective.
Also, if the smallest index at which $q$ and $\hat{q}$ differ is $k$,
then the smallest index at which $f$ and $\hat{f}$ differ is also $k$;
additionally, $f(k) - \hat{f}(k) = q(k) - \hat{q}(k)$.
\end{remark}

We therefore have the following bijection.
\begin{lemma}
\label{bijection}
$F(\Q_n) = \F_n$ and $Q(\F_n) = \Q_n$.
\end{lemma}
\begin{proof}
The map $F : \Q_n \to \F_n$ is surjective because $F(\Q_n) = \F_n$ by the definition of
$\F_n$, and injective, by Lemma~\ref{injective}.
\end{proof}
Hence, there are $n!$ elements in both $\F_n$ and $\Q_n$. We use this
bijection many times in the rest of the paper. In particular, we use it to
`argue in reverse', for instance, proving facts about $q\in\Q$, which may
be easier, to deduce facts about $f\in\F$.

\subsection{Bounds on $f\in\F$}
We first give bounds on $q' = Q'(q)$ for any $q\in\Q$.
\begin{lemma}
\label{q'bounds}
Let $q\in\Q$. Then $q' = Q'(q)$ obeys $q'(1) = q'(2) = 1$ and, for $n\geq 3$, $1\leq q'(n) \leq n-2$.
\end{lemma}
\begin{proof}
The lower bound for $q'(n)$ is easily seen to be 1. For the upper bound,
since $q\in\Q$, $1\leq q(n-1)\leq n-1$, and so $1\leq n-q(n-1)\leq n-1$, both for $n>1$.
Also because $q\in\Q$, we would expect the maximum value of $q'(n)$ to occur for the
maximum allowed value of $n-q(n-1) = n-1$, but this is impossible since
$n-q(n-1) = n-1$ requires $q(n-1) = 1$ and so $q'(n) = 1$ as well. Therefore,
to maximise $q'(n)$, we must have $q(n-1) =2$, so that $n-q(n-1) = n-2$,
whereupon $q'(n)\leq q(n-2)\leq n-2$. 
\end{proof}

A bound on the terms of $f\in\F$ is given in the following.
\begin{lemma}
\label{fbounds}
If $f\in\F$ then $f(1) = 0$, $f(2)\in\{0,1\}$ and for $n>2$, $3-n\leq f(n)\leq n-1$.
\end{lemma}
\begin{proof}
First, $f(1) = 0$ by definition.
Second, since $q(1)=1$, we have $f(2) = q(2) - q(2 - q(1)) = q(2) - q(1) = q(2)-1\in\{0, 1\}$, using~\eqref{qb}.

For $n\geq 3$, we have, from Lemma~\ref{q'bounds}, that $1\leq q'(n)\leq n-2$
and, since $1\leq q(n)\leq n$ and $f(n) = q(n) - q'(n)$, we have $3-n\leq f(n)\leq n-1$.
\end{proof}

In the cases $n=1,2$, it is clear that the bounds in Lemma~\ref{fbounds} are sharp. In
fact they are sharp for all $n$ as we now show.
\begin{lemma}
The bounds given in Lemma~\ref{fbounds} for $n>2$ are sharp. 

For the upper bound, define $f_u(i) = i-1$ for $i=1,\ldots, n$. Then 
$f_u\in\F_n$ and $f_u(n)$ equals the upper bound.

For the lower bound, define the first $n$ terms of sequence $f_l$ by
$$f_l :=(\underbrace{0,\ldots, 0,}_{n-3\mbox{ \scriptsize{zeroes}}} n-3, 1, 3-n).$$
Then $f_l\in\F_n$ and $f_l(n) = 3-n$, which equals the lower bound.
\end{lemma}
\begin{proof}
The upper bound is proved in Lemma~3.2 in~\cite{dilh}. 

Turning to the lower bound, by definition, $f_l(n) = 3-n$, but it remains to prove that $f_l\in\F_n$.
This can be done by direct calculation of $q_l = Q(f_l)$ and checking that
each term in $q_l$ obeys \eqref{qb}. Clearly, $q_l(i) = 1$ for $i =
1,\ldots, n-3$. Furthermore,
$$q_l(n-2) = q_l(n-2-q_l(n-3)) + f_l(n-2) = q_l(n-3) + n-3 = 1 + n-3 = n-2;$$
$$q_l(n-1) = q_l(n-1-q_l(n-2)) + f_l(n-1) = q_l(1) + 1 = 2; \mbox{ and}$$
$$q_l(n) = q_l(n-q_l(n-1)) + f_l(n) = q_l(n-2) + 3-n = n-2+3-n = 1,$$
and so $q_l\in\Q_n$; thus, $f_l\in\F_n$.
\end{proof}

\subsection{When slow $q$ implies slow $F(q)$}

In~\cite{dilh} we gave a few isolated examples of slow $f$ for which $Q(f)$
was also slow. The following lemma adds to these examples.
\begin{lemma}
\label{slow_fq}
Let $\alpha\in \left[\left. \frac{1}{2}, 1\right)\right.$
and define $q(n) = 1 + \lfloor{\alpha n}\rfloor$. 
Then $q(1) = 1$ and $q$ is slow.
Furthermore, $F(q)$ is also slow.
\end{lemma}
\begin{proof}
That $q(1) = 1$ and $q$ is slow is obvious from $\alpha\in \left[\left. \frac{1}{2}, 1\right)\right.$.
To prove that $f = F(q)$ is slow requires us to show that
$f(n+1)-f(n)\in\{0,1\}$ for all $n\in\mathbb N$. That is, that
$$D(n) := \lfloor{\alpha(n+1)}\rfloor - \lfloor{\alpha n - \alpha\lfloor\alpha n\rfloor}\rfloor
-\lfloor{\alpha n}\rfloor + \lfloor{\alpha(n-1) - \alpha\lfloor\alpha(n-1)\rfloor}\rfloor
\in\{0,1\}.$$
First, consider the case $\alpha = 1/2$, which is shown in Lemma~4.1 in~\cite{dilh}
to give $f(n) = \lfloor{(n+2)/4}\rfloor$, clearly a slow sequence.
Hence, we assume that $\alpha\in(1/2, 1)$. We then
write $K = \lfloor{\alpha(n-1)}\rfloor$ and distinguish three cases:
\begin{enumerate}
\item $\lfloor{\alpha n}\rfloor = K$, $\lfloor{\alpha(n+1)}\rfloor = K+1$ 
\item $\lfloor{\alpha n}\rfloor = K+1$, $\lfloor{\alpha(n+1)}\rfloor = K+1$ 
\item $\lfloor{\alpha n}\rfloor = K+1$, $\lfloor{\alpha(n+1)}\rfloor = K+2$.
\end{enumerate}
The bounds on $\alpha$ guarantee that these three cases cover all possibilities.

\textbf{Case (1)} Here,
$$D(n) = K+1 - \lfloor{\alpha(n-K)}\rfloor -K + \lfloor{\alpha(n-1-K)}\rfloor
= 1-\lfloor{\alpha(n-K)}\rfloor + \lfloor{\alpha(n-K)-\alpha}\rfloor.$$
Now, either $\lfloor{\alpha(n-K)}\rfloor = \lfloor{\alpha(n-K)-\alpha}\rfloor$ or
$\lfloor{\alpha(n-K)}\rfloor = \lfloor{\alpha(n-K)-\alpha}\rfloor+1$.
Hence, $D(n)\in\{0, 1\}$.

\textbf{Case (2)} We have
$$D(n) = K+1 - \lfloor{\alpha(n-K-1)}\rfloor -K-1 + \lfloor{\alpha(n-1-K)}\rfloor = 0.$$

\textbf{Case (3)} In this case,
$$D(n) = K+2 - \lfloor{\alpha(n-K-1)}\rfloor -K-1 + \lfloor{\alpha(n-1-K)}\rfloor = 1.$$

The union of the three cases gives $D(n)\in\{0, 1\}$ for all $n$ and therefore $f$ is slow.
\end{proof}
The converse of Lemma~\ref{slow_fq} does not apply: that is, slow $f$
does not imply slow $Q(f)$.

\begin{remark}
We can, however, use the bijection in Lemma~\ref{bijection} to reverse the argument
in Lemma~\ref{slow_fq} to generate a set of slow sequences
$f$ for which $Q(f)$ is also slow. We have:

Let $\alpha\in[1/2, 1)$ and define $f(n) := \lfloor{\alpha n}\rfloor -
\lfloor{\alpha(n-1) - \alpha\lfloor{\alpha(n-1)}\rfloor}\rfloor$. Then $f$
and $q = Q(f)$ are slow and $q(n) = 1 +\lfloor\alpha n\rfloor.$
\end{remark}

\subsection{A family of monotonic $f$ giving monotonic $Q(f)$}
In a similar vein to the previous lemma, we construct a family of
monotonically increasing, but not necessarily slow sequences $f$, for which
$Q(f)$ is also monotonically increasing.

\begin{lemma}
\label{monot_not_slow}
Fix $\delta\in\mathbb N$ and define $f(n) = \delta\lfloor{(n-1)/\delta}\rfloor$.
Then $Q(f) = 1+\delta\lfloor{(n-1)/\delta}\rfloor$ and both $f$ and $q$ are
monotonic but, for $\delta > 1$, not slow.
\end{lemma}
\begin{proof}
This is by direct calculation. We start with $q(n) = 1+\delta\lfloor{(n-1)/\delta}\rfloor$
and deduce the given expression for $f$. We have 
$n-q(n-1) = n-1 - \delta\lfloor{(n-2)/\delta}\rfloor$ and so
$q'(n) = 1 + \delta\lfloor{(n-2)/\delta - \lfloor{(n-2)/\delta}\rfloor}\rfloor$.
Now, for positive $x\in\mathbb R$, we have $x = \lfloor{x}\rfloor + \{x\}$,
where $\{x\}$ is the fractional part of $x$. Hence, $x-\lfloor x\rfloor\in
[0, 1)$, so $q'(n) = 1$ for all $n > 1$, from which we deduce that 
$f(n) = \delta\lfloor{(n-1)/\delta}\rfloor$ as required.
Furthermore, as can easily be checked, we have $q(1) = 1$, $q'(1) = 1$ and $f(1) = 0$
for all $\delta\in\mathbb N$.
\end{proof}

\section{Structure of $\F$}\label{str_F}

The set $\Q$ is easily defined --- in a single sentence --- in such a way
that it is trivial to determine if any given sequence is or is not a
member of $\Q$: just check whether~\eqref{qb} holds for all the terms given.
The set $\F$ is a different matter. So far, all we know is that, given a
sequence $f$, whose terms obey Lemma~\ref{fbounds}, we must use~\eqref{rec} to 
find $q = Q(f)$ and only then can we check whether $q\in\Q$. If it is, then $f\in\F$. 

A natural question would be: is it possible to tell directly whether a given sequence 
$f\in\F$, that is, without going via $Q(f)$? We presently have no neat answer to this,
one reason being that we first need to better understand the structure of $\F$.

To this end, we start with some simple observations: (i) $\F_n$ contains exactly $n!$ sequences of
length $n$ (because the sequences in $\F$ are in one-to-one correspondence
with those in $\Q$), and (ii) the bounds on the terms of a sequence in $\F$ are
given by Lemma~\ref{fbounds}. Item (ii) here already indicates that the
elements of $\F$ have a more complicated structure than those of $\Q$. In $\Q$, the
bounds alone define the sequence, whereas this is not true for $\F$.
There are `gaps' in $\F$, that is, there are sequences that obey 
Lemma~\ref{fbounds} but which are not in $\F$. Let us refer to
the sequences of length $n$ that obey the bounds given by Lemma~\ref{fbounds}
as `potential members of $\F$'. The number of potential members of $\F$ of
length $n$, $\varphi(n)$, is $\frac{(2n-3)!}{2^{n-3}(n-2)!}$, for $n>2$, which, by Stirling's 
approximation, exceeds $n!$ by a factor of about $2^n/\sqrt{\pi n^3}$ for large $n$.

\begin{figure}
\input{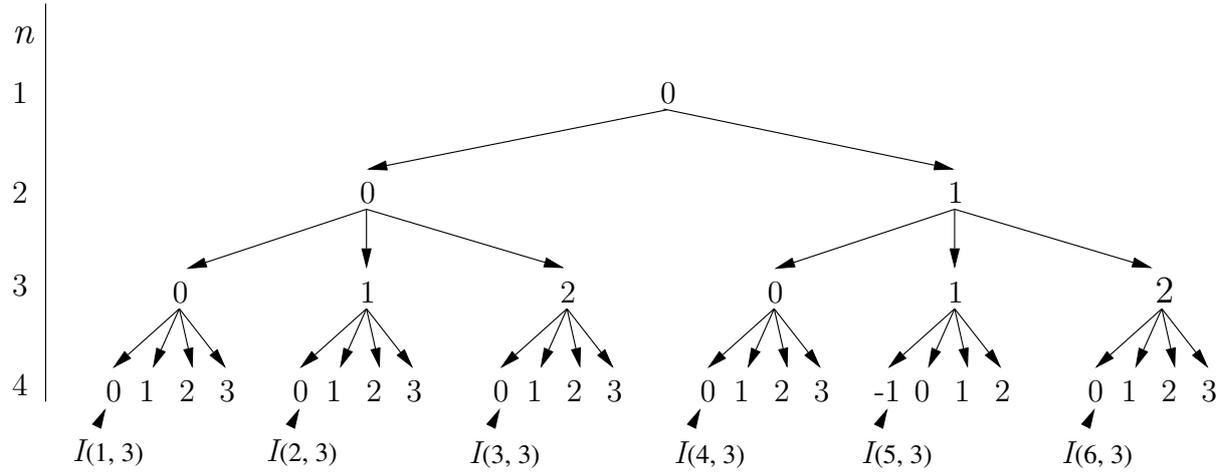}
\caption{The node values in the directed tree structure of $\F_4$. All $4! = 24$ sequences in $\F_4$ are contained
in this diagram. We can see immediately that, for instance,
$(0,1,1,-1)\in\F_4$ but $(0,1,2,-1)\notin\F_4$.}
\label{F_4}
\end{figure}

The set $\F$ can be visualised as a directed tree and the structure of
$\F_4$ is shown in Figure~\ref{F_4}.
Each node in $\F$ has an associated integer value.
The sequence $f$ corresponding to a given node in the $n$-th row of the tree is the sequence of 
the values of the nodes on the path between the root node and the given node.

We define $F_{1,1}$, with value 0, as the root node.
A given node of $\F$ is referred to as $F_{t,n}$ in which $n\in\mathbb N$ relates to
the depth of the node --- specifically, it is $1+$ the number of edges between the root
node and $F_{t, n}$ --- and the integer $t$, $1\leq t\leq n!$ corresponds to
the position of the node within the $n$-th row. At depth $n$, we write the value of
the nodes $F_{in+1, n}$, $0\leq i\leq (n-1)!-1$, as $I(i, n-1)$. The values
of $I(i,3), i = 1,\ldots, 6$, are given in Figure~\ref{F_4}. Note that some sequences 
in $\F$ can contain negative integers: as we see from Figure~\ref{F_4}, $F_{17, 4} = I(5,3) = -1$.

\subsection{Potential members of $\F$}

Consider sequences $f$, each of whose terms obeys the bounds in Lemma~\ref{fbounds},
as potential members of $\F$. Theorem~\ref{fneg} below shows that we can eliminate a significant proportion of such
sequences --- namely, those in which $f(n) < 0$ for all $n$ greater than a given value $n_0 \geq 3$.
The lower bound on $n_0$ just comes from Lemma~\ref{fbounds}, which shows that $f(n)$ cannot be
negative for $1\leq n\leq 3$.

Theorem~\ref{fneg} not only shows that $f\notin\F$ if $f(n) < 0$ for all $n> n_0$,
but also gives an upper bound, $D$, for $d$, the index at which $q := Q(f)$ dies:
that is, we find $D$ such that $q(d) < 1$ --- $q(d) > d$ turns out to be
impossible --- for some $d\leq D$. We give an 
expression for $D$ in terms of $n_0$ and $B := \max_{1\leq n \leq n_0} q(n)$.

Recall that $q'(n) = q(n-q(n-1))$, with $q'(1) := 1$. In order to find the upper bound, $D$,
for the index, $d$, at which $q$ dies,
we need to assume that $q$ exists for as long as possible --- that is,
until the continued existence of $q$ leads to a contradiction --- and in
what follows, we use the abbreviation $E(n)$ to mean: $\left(q(i)\right)_{1\leq i\leq n}\in\Q_n$.

We start by setting up the problem. Let $n_0\geq 3$ and let
$\tilde{f} := (\tilde{f}(n))_{1\leq n \leq n_0}\in \F_{n_0}$. Define
$\tilde{q} := Q(\tilde{f})$. By Lemma~\ref{bijection}, $\tilde{q}\in\Q_{n_0}$.
Define $B := \max_{1\leq n\leq n_0} \tilde{q}(n)$; since
$\tilde{q}\in\Q_{n_0}$, we have $1\leq B \leq n_0$. Let
$f$ be an extension of $\tilde{f}$, where $f(n) = \tilde{f}(n)$ for $1\leq n
\leq n_0$, with $f(n) \leq -1$ for $n > n_0$. Finally, let $q := Q(f)$.

We now prove a preliminary lemma concerning the boundedness of $q$. The
assumptions and definitions set out in the preceding paragraph are all
assumed to hold.
\begin{lemma}
\label{neg_f_q_bound}
Fix $n_0\geq 3$. If $q(n)\leq B$ for $1\leq n\leq n_0$ and $f(n)\leq -1$ for
$n>n_0$, then $q(n)\leq B-1$ for $n > n_0$, as long as $E(n-1)$ holds.
\end{lemma}
\begin{proof}
We use strong induction, and the base step comes from bounding $q(n_0+1) =
q(n_0+1-q(n_0))+f(n_0+1)$. Now, $1\leq q(n_0)\leq B$, and, since $B\leq
n_0$, we have $1\leq n_0+1-q(n_0)\leq n_0$. Hence, $q'(n_0+1) \leq B$. 
Therefore, since $f(n_0+1)\leq -1$, $q(n_0+1) \leq B-1$.

If $B=1$, there is nothing else to prove. If $B>1$,
along similar lines we can now prove the induction hypothesis, which is:
for $j\geq 1$, $q(n) \leq B-1$ for $n_0+ 1 \le n \le n_0+j$ implies that $q(n_0+j+1)\leq B-1$. Note that
$q(n_0+j+1) = q(n_0+j+1-q(n_0+j))+ f(n_0+j+1)$. We have that $n_0+j+1-(B-1)\leq n_0+j+1 - q(n_0+j) \leq n_0+j$, 
using the bounds on $q(n_0+j)$. Hence, since
$n_0+j+2-B$ may or may not be greater than $n_0$, the most we can say,
without further analysis, is
that $q'(n_0+j+1)\leq B$, and so, since $f(n_0+j+1)\leq -1$, we have
$q(n_0+j+1)\leq B-1$.
\end{proof}
This simple lemma gives insight into the argument that we will use to prove Theorem~\ref{fneg}.
Under the conditions that $q(n)\leq B$ for $1\leq
n\leq n_0$ and $f(n) \leq -1$ for $n > n_0$, we have that, as $n$
increases beyond $n_0$, $n-q(n-1)\geq n-B+1$ must also increase. Hence, there must
exist an integer $n$ for which $n-q(n-1) > n_0$, and this forces $q(n-q(n-1))\leq B-1$, by
Lemma~\ref{neg_f_q_bound}. Therefore, from this point on, $f(n)\leq -1$ gives $q(n)\leq B-2$.

Repeating the above argument, we will find an ever-decreasing sequence of upper 
bounds on $q(n)$ as $n$ increases. Specifically, for given
$n_0$ and $B$, there must exist a sequence of integers
$n_k$, with $k \in \{ 1:B\}$, obeying $3\leq n_0 < n_1< \ldots < n_k < \ldots < n_B$, such
that, given that $1\leq q(n)\leq B$ for $1\leq n\leq n_0$, then,
for $k \in \{1:B\}$ and $n \geq n_k$, either $E(n-1)$ does not hold or $q(n)\leq B-k$.
We derive  an expression for the $n_k$ in Theorem~\ref{fneg}, and this then allows us to compute $D$.

From the preceding argument, it is clear that, first, $q\notin \Q$ and,
second, $q$ dies at some $d$ because $q(d) < 1$ --- and not because $q(d) > d$.
Indeed, on this last point, at worst, if the sequence $q$ is well-defined up to $n=n_B-1$, 
then $q(n_B) \le 0$ and hence $d\le n_B$;
on the other hand, by Lemma~\ref{neg_f_q_bound}, $q(n)\leq B-1\leq n_0-1$ for $n>n_0$, so we
cannot have $q(d)>d$ for any $d>n_0$.

It is important to note that in proving Lemma~\ref{neg_f_q_bound}, we
assumed only that $f(n)\leq -1$. The argument does not change at all if, for some
values of $n$, $f(n)<-1$. The upper bound on $q(n)$ would only change if, for
instance, for all $n>n_0$, $f(n)\leq -c$ for some $c > 1$. In that case, the bound would be
modified to $q(n)\leq B-c$ for $n > n_0$. In order to compute $D$, however,
we wish $q$ to survive for as long as possible, and therefore we use $f(n) = -1$ for
$n>n_0$ when proving the following.

\begin{theorem}
\label{fneg}
Fix $n_0\geq 3$ and let $\tilde{f} := (\tilde{f}(n))_{1\leq n \leq n_0}\in \F_{n_0}$.
Define $\tilde{q} := Q(\tilde{f})$ and $B := \max_{1\leq n\leq n_0} \tilde{q}(n)$. Let
$f$ be an extension of $\tilde{f}$, where $f(n) = \tilde{f}(n)$ for $1\leq n \leq n_0$,
with $f(n)\leq -1$ for $n > n_0$. Define $q := Q(f)$.
Then
\begin{enumerate}[label=(\roman*)]
\item There exists a finite sequence of integers $n_0 < n_1 <  \ldots < n_k < \ldots < n_B$, 
such that, for $k=1,\ldots,B$, we have
\begin{subequations} \label{ind_hyp}
\begin{align}
& n_k = n_{k-1}+B-(k-1),\mbox{ and}
\label{ind_hypa} \\
& q(n)\leq B-k \quad \forall n \ge  n_k \hbox{ as long as } E(n-1) \hbox{ holds.}
\label{ind_hypb}
\end{align}
\end{subequations}
\item There is an integer $d$ such that $q(d) < 1$, so that $q$ dies at $d$, where $d\le D$, with
$$ D = D(n_0,B) :=  n_0+\frac{1}{2}B(B+1) . $$
\end{enumerate}
\end{theorem}
\begin{proof}
If $B=1$, the result follows directly from Lemma \ref{neg_f_q_bound}, with $d=n_1=n_0+1$. If $B>1$, we proceed as follows.

For (i), we start by considering $k=1$, and this will be the base case for (strong) induction.
By setting $n_1=n_0+B > n_0$, the inequality $1\le q(n) \le B-1$ for $n \ge n_1$ follows from 
Lemma \ref{neg_f_q_bound}.
Next, consider $k\in\{2:B\}$ and assume~\eqref{ind_hypa}
and~\eqref{ind_hypb} for all $k'=1,\ldots,k-1$.  For $n\ge n_{k}$, we have
$n-1 \ge n_k - 1 = n_{k-1} + B - (k-1) - 1 = n_{k-1} + B - k \ge n_{k-1}$ and hence, by the 
induction hypothesis~\eqref{ind_hyp}, $q(n-1) \le B-(k-1)$.
This implies that $n-q(n-1) \ge n - B +(k-1) \ge n_k - B + (k - 1) = n_{k-1}$, so 
that $q(n-q(n-1)) \le B- (k-1)$, which in turn gives $q(n)=q(n-q(n-1)) +f(n) \le B- (k-1) - 1 \le B-k$.

For (ii), the recursion formula \eqref{ind_hypa} for $n_k$ is easily solved to give 
$n_k = n_0 +  kB - k(k-1)/2$, so that $n_B=n_0+B(B+1)/2$. Since $q(n_B) \le 0$,
the smallest $n$ for which $q(n)< 1$ is no larger than $n^* = n_0 + B(B+1)/2$.
Therefore, $D = n^*$ is a bound on the value $d$ at which the sequence $q$ dies.
\end{proof}

\begin{remark}
Theorem~\ref{fneg} allows us to eliminate a significant number (i.e., a number that
grows factorially rather than, say, exponentially) of potential members of $\F$,
being those all of whose terms are negative for some $n_0>3$. Considering
sequences $f$ of length $n>3$, the number which consist of only negative terms from the
fourth term onwards is $(n-3)!$. The proportion of the $\varphi(n)$ potential
sequences eliminated by Theorem~\ref{fneg} is therefore
$$\frac{(n-3)!}{\varphi(n)}
\sim\frac{1}{2^n}\sqrt{\frac{\pi}{n^3}}\left(1 + O\left(\frac{1}{n}\right)\right).$$
\end{remark}

It is possible to do an exhaustive search for $4\leq n_0\leq 12$, in order to
find, for each $n_0$ the value of $B$, $B_{\mathrm{comp}}$, which leads to the largest $d$.
We proceed as follows: (1) generate all $n_0!$ sequences $f\in\F_{n_0}$;
(2) extend each with $f(n) = -1$ for $n>n_0$, noting where each sequence dies, which happens at $n=d$;
(3) note the maximal value of $d$ and count how many of the $n_0!$ sequences die at this maximal $d$.
The count in (3) is $N_D$ and the maximal value of $d$ is $D_{\mathrm{comp}}$.
The results are given in Table~\ref{tabD}.
All $N_D$ sequences for each $n_0$ have $B = n_0-1$,
suggesting that this is always the case. Certainly, it is clear that $B = n_0$
can only occur because $\tilde{q}(n_0) = n_0$, so that $q(n_0+1) =
q(n_0+1-n_0)+f(n_0+1) = 1 + f(n_0+1)\leq 0$. Hence, $B = n_0$ implies that $d = n_0+1$,
and therefore for maximal $d$ we must have $B\leq n_0-1$.
\begin{table}[h!]
\centering
\begin{tabular}{|l|c|c|c|c|c|c|c|c|c|}\hline
$n_0$    &      4   &  5  &  6  &  7  &  8  &  9  &  10  &  11  &  12\\ \hline
$B_{\mathrm{comp}}$ & 3  &  4  &  5  &  6  &  7  &  8  &   9  &  10  &  11\\ \hline
$D_{\mathrm{comp}}$ & 7  &  9  &  14 &  22 &  27 &  37 &  46  &  56  &  67\\ \hline
$D(n_0,B)$    & 10 &  15 &  21 &  28 &  36 &  45 &  55  &  66  &  78\\ \hline
$N_{D} $    & 2  &  14 &  6  &  40 &  40 & 280 & 280  & 9880 &  34520\\ \hline
\end{tabular}
\vspace{0.5ex}
\caption{\label{tabD} Comparison of $D$ computed by exhaustive search, and as
given in Theorem~\ref{fneg}.}
\end{table}

\subsection{Some subsets of $\F$}
\label{subsets}
It was proved in~\cite{dilh} that for all slow sequences $f$ with $f(0)=0$, $Q(f)$
exists. The set of all such slow sequences is therefore an
easily generated subset of $\F_n$ that contains exactly $2^{n-1}$ sequences.
Exponential growth may seem impressive, but in the current context
this is a vanishingly small proportion of the $n!$ sequences in $\F_n$. Can we do better?

In order to compare different subsets of $\F_n$, we introduce the function
$C(n)$, which gives the number of sequences of length $n$ in a given subset
of $\F_n$ as $n$ varies --- or an approximation to this number. 
For the set of slow sequences, for example, we have $C(n) = 2^{n-1}$.
As the following examples will show, we are mainly interested in whether $C(n)$ grows exponentially or faster
than exponentially, and $C(n)$ needs to allow us to distinguish between these two cases.

We start by giving a construction of a set of sequences that are, intuitively, far from slow.
\begin{lemma}\label{pwr2}
Let $f(n)\in\{0, n-1\}$ for $n\in\mathbb N$ and let $q = Q(f)$. Then
$q(n) = 1+f(n)\in\{1, n\}$ and so $q\in\Q$, and therefore $f\in\F$.
\end{lemma}
\begin{proof}
We use induction. By definition, we have $f(1) = 0$ and $q(1) = 1$, so
consider $q(2) = q(2 - q(1)) + f(2) = 1+f(2)$ as the base case.

Now let $k>2$ and assume that $\left(q(i)\right)_{1\leq i\leq k}= \left(1 +f(i)\right)_{1\leq i\leq k}$.
Then
$$q(k+1) = q(k+1 - q(k)) + f(k+1) = \begin{cases}
q(1) + f(k+1) & \mbox{if}\qquad q(k) = k\\
q(k) + f(k+1) & \mbox{if}\qquad q(k) = 1.
\end{cases}$$
Hence, in either case, $q(k+1) = 1 + f(k+1)\in\{1, k+1\}$, so
$q\in\Q$ and therefore $f\in\F$.
\end{proof}
This lemma, which turns out to be a special case of Theorem~\ref{fastest},
gives examples of sequences $f\in\F$ in which the difference between
successive terms, $f(n+1)-f(n)$, can be made as large as we like by
choosing sufficiently large $n$. Since the choice for $f(n)$ can be made
freely from $\{0, n-1\}$ for each $n$ except $n=1$, we have $C(n) = 2^{n-1}$.


Returning to the question of subsets of $\F_n$ whose size grows faster than
$2^n$, Lemma~\ref{fbounds} leads to the example in Lemma~\ref{Fsubset3}
below. Before stating and proving the lemma, we discuss the
following question. Given $n>1$, suppose we have a sequence $\tilde f\in\F_n$; how
do we extend $\tilde f$ by one term so that the resulting sequence, $f$, is in $\F_{n+1}$?
To answer this, let $\tilde q := Q(\tilde f)$ and $q := Q(f)$, so that $q(n+1) = q(n+1-q(n))+ f(n+1)$
is well-defined. We must have $q(n+1)\in\{1:n+1\}$ because $q\in\Q_{n+1}$.
Now, considering \textit{all} valid extensions $q$ of $\tilde q$, we can choose $q(n+1)$ freely 
from $\{1:n+1\}$, and the choice made does not affect the value of
$q'(n+1) := q(n+1- q(n))$, which we can consider to be fixed. Hence,
$f(n+1)$ lies in a difference set, specifically,
\begin{equation}
\label{fn1_set}
f(n+1)\in\{1:n+1\} - \{q(n+1-q(n))\}.
\end{equation}
In words, $f(n+1)$ lies in the set of $n+1$ consecutive integers $\{1- q(n+1-q(n)): n+1-q(n+1- q(n))\}$. 

\begin{lemma}
\label{Fsubset3}
For all sequences $f$ that obey
$$f(n) \in \begin{cases}
	\{0\} & n = 1\\
	\{0, 1\} & n = 2\\
	\{0, 1, 2\} & n \geq 3
\end{cases}$$
$Q(f)$ exists.
\end{lemma}
\begin{proof}
It is clear that all $f\in\F$ must obey $f(1) = 0$ and either $f(2) = 0$ or $f(2) = 1$.
For $f(n)$, $n\geq 3$, Lemma~\ref{fbounds} gives $3-n\leq f(n)\leq n-1$,
and, as implied by~\eqref{fn1_set}, there are $n$ allowed consecutive integer values of $f(n)$.
Hence, in order to respect both these constraints, for $n\geq 3$, $f(n)$ must lie in the set
$\cap_{i=0}^{n-3} \{3-n+i: 2+i\} = \{0:2\}$, which completes the proof.
\end{proof}
By inspection, for $f$ as in Lemma~\ref{Fsubset3} we have $C(n) = 2\cdot 3^{n-2}$ for $n\geq 3$.
\begin{corollary}[Extension of a sequence in $\F_{n_0}$] 
Let $f_{n_0}\in\F_{n_0}$ be a sequence of $n_0\geq 3$ terms.
Let the infinite sequence $f$ be defined by $f(n) = f_{n_0}(n)$, $1\leq n \leq n_0$,
with $f(n)\in\{0:2\}$ for $n > n_0$. Then $f\in\F$.
\end{corollary}
\begin{proof}
This is by induction. Fix $k > n_0$ and assume that $f_k := \left(f(i)\right)_{1\leq i\leq k}\in\F_k$.
Since $f_k\in\F_k$, $Q(f_k)\in\Q_k$.  Therefore $q(k+1) = q(k+1-q(k)) + f(k+1)$
is well-defined, and $1\leq q(k+1-q(k))\leq k-1$, because either $q(k)=1$, whereupon $q(k+1-q(k))=q(k)=1 < k-1$;
or $q(k)>1$, so that $q(k+1-q(k)) \le \max_{1 \le i \le k-1} q(i) \le k-1$.
Finally, since $f(k+1)\in\{0: 2\}$, $1\leq q(k+1)\leq k+1$.
\end{proof}

We can do better than $C(n) = 2\cdot 3^{n-2}$, however, as the following lemma shows. The proof of
this, and of Theorem~\ref{fastest}, follow the same line of argument, which is: a set of sequences $q\in\Q$ are
proposed, these being such that there is a simple expression for $q' = Q'(q)$,
and hence, by subtraction, for $F(q)$. In the light of the
bijection~\ref{bijection}, we are in effect exploiting the fact that we can
prove things either starting with $q$ and deducing properties of the 
corresponding $f$, or vice versa.

\begin{lemma}
\label{l_strip}
Fix $\ell\geq 3$ and consider the sequences $q$ defined in the row labelled
`$q(n)$' in Table~\ref{qband}.
\begin{table}[h!]
\centering
\begin{tabular}{|l|c|cc|l|c|cccc|}\hline
$n$         & $1$ & $2$ & $3$ & $4$\ldots $\ell+1$ & $\ell+2$       & $\ell+3$       &\ldots   & $\ell+k$            & \ldots\\ \hline
\rb $q(n)$  & $1$ & $2$ & $2$ & $2\ldots 2$        & $\{3:\ell+1\}$ & $\{3:\ell+2\}$ &\ldots   & $\{k:\ell+k-1\}$    & \ldots\\
$q'(n)$     & $1$ & $1$ & $1$ & $2\ldots 2$        & $2$            & $2$            &\ldots   & $2$                 & \ldots\\
$f(n)$      & $0$ & $1$ & $1$ & $0\ldots 0$        & $\{1:\ell-1\}$ & $\{1:\ell\}$   &\ldots   & $\{k-2:\ell+k-3\}$  & \ldots\\ \hline
\end{tabular}
\vspace{0.5ex}
\caption{\label{qband} The sequences $q$, $q'$ and $f$ used in Lemma~\ref{qband}.}
\end{table}
Then $q\in\Q$ and $f = F(q)\in\F$.
\end{lemma}
\begin{proof}
Note that (a) for $n\geq\ell+2$ the indicated choices for $q(n)$ can be made
freely; and (b) for $n = \ell+2$, we could also have $q(n) = 2$, but that would be
equivalent to increasing $\ell$ by 1, so we exclude this case. 

It is obvious from the definition that $q\in\Q$.

The sequences $q$ here are all contrived to give $q'(n) = 2$ for $n\geq 4$. This
is because $n-q(n-1) \in \{2:\ell\}$ when $n=\ell+3$, and 
$n-q(n-1) \in \{2:\ell+1\}$ for $n>\ell+3$.
Hence, $f(n) = q(n) - 2$ for $n\geq\ell +2$, with $f(n)$ for $n<\ell+2$
being $0$ or $1$, as given in Table~\ref{qband}.
\end{proof}

The set of the first $n$ terms of the sequences $f$ generated in this way has
$$C(n) = (\ell-1)\ell^{n-\ell-2}$$
members for $n\geq\ell +2$ --- still exponential growth,
but here we can make $\ell$ as large as we please.

Thus, the original question should perhaps be changed to: can we generate a subset of $\F_n$
for which $C(n)$ grows faster than exponentially with $n$? The best result
we have achieved in this direction has `scaled factorial' growth, that is, $C(n)\approx (\alpha n)!$ for $\alpha\in(0, 1)$.
Such a growth rate always exceeds exponential growth: for any $\alpha,
\beta > 0$
$$\frac{(\alpha n)!}{\beta^n}\sim (\gamma n)^{\alpha n}
\qquad\mbox{where}\qquad \gamma = \alpha\beta^{-\frac{1}{\alpha}}\mathrm{e}^{-1} > 0,$$
which tends to infinity with $n$.

Theorem~\ref{fastest} below gives the construction of a subset of $\F_n$
with scaled factorial growth rate. The idea behind its proof is to construct sequences $q$
such that $Q'(q) = \left(q'(n)\right)_{n\in\mathbb N} = 1$ for all $n$.

\begin{theorem}
\label{fastest}
Let $\Y$ be the set of sequences $\left(y(n)\right)_{n\in\mathbb N}$ with
$y(1) = 1$ and $y(n)\in\{0, 1\}$ for $n>1$. With each sequence $y\in \Y$
associate an infinite sequence of sets, $S(y) := \left(s(n)\right)_{n\in\mathbb N}$, by
\begin{equation}
\label{Sdef}
s(n) :=
\begin{cases}
\{1\} & \mbox{if $y(n) = 1$}\\
\{1\} \cup \left(\{n+1\} - \{j | y(j) = 1, 1\leq j\leq n\}\right) & \mbox{if $y(n) = 0$.}
\end{cases}
\end{equation}
Let sequences $q$ be formed by choosing each member $q(n)$ from $s(n)$, independently, for each
$n$. Then all such $q$ are members of $\Q$ and, defining
$q'(n) := q(n-q(n-1))$ for $n\in\mathbb N$, with $q'(1) = 1$, we have $q'(n) = 1$ for all $n$.
\end{theorem}
\begin{proof}
Let $y\in \Y$ and $s = S(y)$ be as defined in equation~\eqref{Sdef}. Let us
write $q\in S(y)$ to mean $q(n)\in s(n)$ for all $n$. To show that, for all
$y\in\Y$, any $q\in S(y)$ is a member of $\Q$, consider the bounds on $q(n)$. First, note that
$\max\{j| y(j) = 1, 1\leq j \leq n\} \leq n$; and that $\min\{j| y(j) = 1, 1\leq j \leq n\} = 1$,
since $y(1) = 1$. Hence, $\max s(n) = \max(\{n+1\} - \{1\}) = n$ and
$\min s(n) = \min(\{n+1\} - \{n\}) = 1$. Thus, $1\leq q(n) \leq n$ for all $n$
and so $q\in\Q$.

Second, to show that $q' = Q'(q) = \left(q'(n)\right)_{n\in\mathbb N} = 1$ for any $q\in S(y)$,
recall that $q'(1) = 1$ by definition, so consider $q'(n+1)$, $n\in\mathbb N$. Now, since
$1$ is always a member of $S(y)$, if we choose $q(n) = 1$, we have $q'(n+1) = q(n+1-q(n)) = q(n) = 1$.
Therefore, we only need to consider the second case in~\eqref{Sdef}, excluding $q(n) = 1$. For this case, define
$$\sigma(n+1) := \{n+1\} - \{j | y(j) = 1, 1\leq j\leq n\}$$
and note that the difference set
$$\tau(n+1) := \{n+1\} - \sigma(n+1) = \{j | y(j) = 1, 1\leq j\leq n\}.$$
Now, by definition, $q'(n+1) \in \{ q(j) | j\in \tau(n+1)\} = \{1\}$, since
$\tau(n+1)$ is exactly the set of indices $j$ for which $y(j) = 1$.
Hence, $q'(n+1)=1$ for $n\in\mathbb N$, and this completes the proof.
\end{proof}

We make three remarks on this result.
\begin{enumerate}[label=(\roman*)]
\item Let $y\in \Y$, $s = S(y)$ and $q\in S(y)$, with each choice for
$q(n)\in s(n)$ being made freely. Letting $q' = Q'(q)$ we always have $q'(n) = 1$.
Furthermore, we have $\left(f(n)\right)_{n\in\mathbb N} = q(n) - 1$ for all
$n$. Hence, all sequences $f$ defined by
\begin{equation}
\label{fdef}
f(n) \in
\begin{cases}
\{0\} & \mbox{if $y(n) = 1$}\\
\{0\} \cup \left(\{n\} - \{j | y(j) = 1, 1\leq j\leq n\}\right) & \mbox{if $y(n) = 0$}
\end{cases}
\end{equation}
are members of $\F$.
\item We now consider $C(n)$ for sequences $f$ obeying~\eqref{fdef}. By
inspection, we can give an explicit formula for $C(n)$ for any $y\in \Y$,
this being
\begin{equation}
\label{Cfastest}
C(n) = \prod_{i=1}^n T(i)\qquad\mbox{where}\qquad T(i) := 1 + (1 -
y(i))\sum_{j=1}^i y(j).
\end{equation}
Provided that $\sum_{i=1}^n y(i)$ is unbounded as $n\to\infty$, $C(n)$ will
be bounded below by a scaled factorial. If, by contrast,
$\sum_{i=1}^n y(i) = K$ for $n$ greater than some integer, $C(n)\sim A (K+1)^n$, with
$A$ a constant --- only exponential growth.
\item Let us consider the explicit example
$$y(n) = \begin{cases} 1 & n\equiv 1 \pmod{m}\\ 0 & \mbox{otherwise,} \end{cases}$$
for fixed $m > 1$, so that $y$ is periodic with period $m$.
Equation~\ref{Cfastest} gives
$$C(n) \approx \left(\left(\tfrac{n+m}{m}\right)!\right)^{m-1},$$
which is exact for $n\equiv 0\pmod{m}$. Stirling's approximation for this gives
$$\ln C(n) =  (m-1)\left[\frac{n}{m} \ln\left(\frac{n}{m}\right) -
\frac{n}{m} +\frac{3}{2} \ln\left(\frac{n}{m}\right) + \frac{1}{2} \ln 2\pi\right] + O\left(\frac{1}{n}\right).$$
In this example, $C(n)$ displays scaled factorial growth: recall that $\ln (\alpha n)! \sim
\alpha n (\ln (\alpha n) - 1)$ and substitute $\alpha = 1/m$.
\end{enumerate}

\subsection{Linear cones}
A linear cone, $\C(a, b)$, is a set of infinite sequences and is defined next.
\begin{definition}[Linear cone]
Let $0\leq a < 1$ and $a < b \leq 1$. An infinite sequence $s$ is said to be a member of the linear 
cone $\C(a, b)$ if, for all $n\in\mathbb N$, 
$$
\lfloor{an}\rfloor\leq s(n)\leq 
\begin{cases}
\lfloor{bn}\rfloor & \mbox{for $b\in (a,1)$}\\
n-1 \qquad & \mbox{for $b = 1$}.
\end{cases}$$
\end{definition}
We now prove two negative results concerning linear cones. We require
the following particular set of infinite sequences.
\begin{definition}[Sequences $f_d$]
Fix even, positive integer $d$ and define the set of sequences $\left(f_d(n)\right)_{n\in\mathbb N}$ by
$$f_d = (\underbrace{0,}_{d-1}\underbrace{1,0,}_{d/2}\;\;
\underbrace{2,0,}_{d/2}\;\;
\underbrace{3,\cstar,}_{d/2}\;\;
\underbrace{4,\cstar,}_{d/2}\;\;
\underbrace{5,\cstar,}_{d/2}\;\;
\underbrace{6,\cstar,}_{d/2}\ldots).$$
Here, (a) $\cstar$ is either 0 or 1, and either choice can be made on each
occurrence of $\cstar$; and (b) the number of occurrences of each subsequence is written below 
the underbrace under that subsequence.
\end{definition}
Since $\cstar$ can be 0 or 1, we have defined a set of sequences rather than a single one.

As can easily be checked, (i) all sequences $f_d$ are contained within $\C(0, 1/d)$,
(ii) for even $n>0$, $f(n) = \lfloor{n/d}\rfloor$, and (iii) for odd $n>0$, $f(n)\in\{0,1\}$.

As an example, for $d = 4$ we have
$$f_4 = \left(0, 0, 0,\;\; 1, 0,\; 1, 0,\;\; 2, 0,\; 2, 0, \;\;
3, \cstar_1,\; 3, \cstar_2,\;\; 4, \cstar_3,\; 4, \cstar_4, \ldots \right),$$
where $\cstar_1,\ldots, \cstar_4\in\{0, 1\}$.

\begin{lemma}
\label{cone0}
For all $b\in(0, 1)$, there are sequences $f\in\C(0, b)$ which are not in $\F$.
\end{lemma}
\begin{proof}
We give a proof by construction. Clearly, if the lemma is true for $b$
it is also true for $b'>b$, since $\C(0,b)\subset \C(0,b')$. The idea is to fix an even integer $d\geq 2$ and
then construct the sequence $q_d := Q(f_d)$, for which there exists a finite
integer $K$, dependent on $d$, defined by $q_d(K)>K$.\footnote{It is
tempting to write that `$K$ is the smallest such integer', but `smallest' is
redundant, since $q_d(K+1)$ is undefined.} At this point, the sequence $q_d$ is said to die
because $q_d(K+1)$ is undefined.
For any even $d\geq 2$, it turns out that such a $K$ can always be found.

For the rest of the proof, we write $q_d$ as $q$ and $f_d$ as $f$.
By direct calculation, the first terms of $f$ and $q$ are easily found to be
\begin{table}[h!]
\centering
\begin{tabular}{|l|lcc|ccccc|cccc|cc|}\hline
$n$ & $1$ & $\ldots$ & $d-1$ & $d$ & $d+1$ & $\ldots$ & $2d-2$ & $2d-1$ & $2d$ & $2d+1$ & $\ldots$ & $3d-1$   & $3d$ & $\ldots$\\ \hline
$f$ & $0$ & $\ldots$ & $0$   & $1$ & $0$   & $\ldots$ & $1$    & $0$    & $2$  & $0$    & $\ldots$ & $0$      & $3$  & $\ldots$\\ \hline
$q$ & $1$ & $\ldots$ & $1$   & $2$ & $1$   & $\ldots$ & $2$    & $1$    & $3$  & $2$    & $\ldots$ & $\ldots$ & $\ldots$ & $\ldots$\\ \hline
\end{tabular}
\vspace{0.5ex}
\caption{\label{nfq} Initial values of $f$ and $q$.}
\label{ifq}
\end{table}

In Table~\ref{nfq}, we have $q(2d) = q(2d - q(2d-1)) + f(2d) = q(2d-1) + f(2d) = 3$ and
$q(2d+1) = q(2d+1 - q(2d)) + f(2d+1) = q(2d-2) + f(2d+1) = 2$. We will use
these two values as the base cases for induction.

We now assume, for some $m\geq 1$, that
\begin{equation} \label{base}
q(2d + 2m) = 3 + 2m + \sum_{j=1}^m \left\lfloor\frac{2j}{d}\right\rfloor
\qquad\mbox{and}\qquad q(2d + 2m + 1) = 2,
\end{equation}
and then proceed by induction. We have, for even-indexed $q$,
\begin{equation*}
\begin{split}
q(2d + 2m + 2) &= q(2d + 2m + 2 - q(2d + 2m +1))+f(2d+2m+2)\\
&= q(2d + 2m)+\left\lfloor{\frac{2d+2(m+1)}{d}}\right\rfloor\\
&= 3 + 2m + \sum_{j=1}^m \left\lfloor\frac{2j}{d}\right\rfloor + \left\lfloor{\frac{2d+2(m+1)}{d}}\right\rfloor\\
&= 3 + 2(m+1) + \sum_{j=1}^{m+1} \left\lfloor\frac{2j}{d}\right\rfloor,
\end{split}
\end{equation*}
using~\eqref{base}.

For odd-indexed $q$, we find
\begin{equation*}
\begin{split}
q(2d + 2m + 3) &= q(2d + 2m + 3 - q(2d + 2m +2))+f(2d+2m+3)\\
&= q\left(2d - 2 - \sum_{j=1}^{m+1} \left\lfloor\frac{2j}{d}\right\rfloor\right)+\cstar,
\end{split}
\end{equation*}
where, as before, $\cstar\in\{0, 1\}$. Defining
$$p(d, m) := 2d-2 - \sum_{j=1}^{m+1} \left\lfloor\frac{2j}{d}\right\rfloor,$$
we immediately have $p(d, m) \leq 2d-2$. Hence, from Table~\ref{nfq}, $q\left(p(d,m)\right)\in\{1, 2\}$,
provided that $q(p(d, m))$ is defined: i.e., as long as $p(d, m) \geq 1$.
Hence, by choosing the appropriate
value for $\cstar$, we can ensure that $q(2d + 2i+1) = 2$ for all $i\in\mathbb N$.

The rest of the proof now follows straightforwardly. Clearly, there is a
finite $m$ such that $p(m, d)\leq 0$, and at this point, sequence $q$ dies.

Therefore, we have proved by construction that there is
a sequence $f\in\C(0, 1/d)$ such that $Q(f)$ does not exist. Since $d$ can be
any even positive integer, we can make $b>0$ in $\C(0, b)$ as small as we like,
and from this the lemma follows.
\end{proof}
In this construction, note how the $\cstar$ terms in $f$ are used to
`steer' the sequence $q$ to ever larger values, such that eventually $q$ dies.

With a little more analysis, we can estimate the integer $K = K(d)$ where the sequence $q_d=Q(f_d)$ dies,
and thereby strengthen Lemma~\ref{cone0}.
We first define $r(id),\; i\geq 1$, noting that $q(id) = r(id)$ for $i$ such that $q(id)$ exists.
Going back to equation~\eqref{base} and letting $2m = kd$ with $k\geq 1$, we set

$$r((k+2)d) = 3 + kd + S(d, k)\qquad\mbox{where}\qquad S(d, k) := \sum_{j=1}^{kd/2} \left\lfloor{\frac{2j}{d}}\right\rfloor.$$
Note also, from Table~\ref{ifq}, that $r(d) = 2$ and $r(2d) = 3$.
For $id>K$, $r(id)$ is as defined, whereas $q(id)$ is undefined.

By inspection, $S(d, k) = k\left(1+ \frac{(k-1)d}{4}\right)$. Defining $\ell = k+2$, we then have
\begin{equation}
\label{qld}
r(\ell d) = \begin{cases}
2, 3 & \mbox{$\ell = 1, 2$ respectively}\\
(\ell+1)\left(1 + \frac{(\ell-2)d}{4}\right) & \mbox{$l\geq 3$.}
\end{cases}
\end{equation}
Now we can see that $r(\ell d) = O(\ell^2)$ and therefore there must exist
a finite integer $K$ at which $q$ dies --- this is just Lemma~\ref{cone0} again.

With $K$ defined by $q_d(K) > K$ we can
increase the lower bound for the linear cone in Lemma~\ref{cone0} from 0 to
$a := 1/(K+\varepsilon)$, for any $\varepsilon>0$, this definition of $a$ guaranteeing that $\lfloor{an}\rfloor = 0$
for $1\leq n\leq K$. We then have the following corollary to Lemma~\ref{cone0}.
\begin{corollary}
For $b\in (0, 1)$, $a\in (0,b)$ sufficiently small and $d\geq 2$ an even integer, there are sequences 
$f\in\C(a, b)$ which are not in $\F$.
\end{corollary}
We can use equation~\eqref{qld} to estimate $K \approx xd$ where $x\in\mathbb R$
satisfies $x d = r(x d)$. We find
$$x(d) := \frac{5d - 4 + \sqrt{(11d-4)(3d-4)}}{2d},$$
giving:
\begin{table}[h!]
\centering
\begin{tabular}{|c|llllllllll|}\hline
$d$                         & $2$  & $4$  & $6$  & $8$  & $10$ & $12$ & $14$ & $16$ & $18$ & $20$  \\ \hline
$\vp\lceil{d x(d)}\rceil$   & $6$  & $17$ & $28$ & $39$ & $50$ & $60$ & $71$ & $82$ & $93$ & $103$ \\ \hline
$K(d)$                      & $8$  & $18$ & $30$ & $40$ & $50$ & $62$ & $72$ & $82$ & $94$ & $104$ \\ \hline
\end{tabular}
\end{table}

This table compares the actual value of $K(d)$ with the approximation $\lceil{d x(d)}\rceil$. 

As an example, consider $d = 6$, for which
$$f = (0, 0, 0, 0, 0, 1, 0, 1, 0, 1, 0, 2, 0, 2, 0, 2, 0, 3, 1, 3, 0, 3, 1, 4, 1, 4, 1, 4, 1, 5)$$
and
$$Q(f) = (1, 1, 1, 1, 1, 2, 1, 2, 1, 2, 1, 3, 2, 5, 2, 7, 2, 10, 2, 13, 2, 16, 2, 20, 2, 24, 2, 28, 2, 32).$$
We have here that $K(6) = 30$ and $q(30) = 32$.

The following lemma gives a related result.
\begin{lemma} \label{cone1}
For any $\varepsilon\in (0,1)$ there are sequences $f\in\C(1-\varepsilon, 1)$ which are not in $\F$.
\end{lemma}
\begin{proof}
Clearly, if the lemma is true for some $\varepsilon \in(0,1)$ then it is also
true for any $\varepsilon' \in(\varepsilon, 1)$. Hence, we need only consider
$\varepsilon\to 0$. Again, the proof is by construction.

Fix positive integer $N$ and define $f(i) = i-1$ for $1\leq i\leq N$,
$f(N+1) = N-1$ and $f(N+2) = N+1$. Let $q = Q(f)$. Then $q(i) = i$ for
$1\leq i\leq N$, $q(N+1) = q(N+1-q(N)) + N-1 = N$ and $q(N+2) =
q(N+2-N)+N+1 = N+3$. Hence, $q$ dies at $N+2$. 

Now, the sequence $f$ as defined is in $\C(1-\varepsilon, 1)$ with
$\varepsilon >\frac{1}{N+1}$.
Hence, we can make $\varepsilon$ as small as we like by choosing a
sufficiently large $N$.
\end{proof}

Thus far, we have two negative results concerning linear cones. It may even be the
case that there are no linear cones which form a subset of $\F$ --- and
condition (3) in Theorem~\ref{thm:2} explicitly forbids them --- but we have no proof.

\section{Two further subsets of $\F$}\label{further_subsets}

In the previous section, many of the results were proved by construction,
by which we mean that an explicit construction was devised for, say, a set of sequences,
and the proof proceeded from there.  By contrast, in this section, we adopt an analytical approach
to prove two further results. Theorem~\ref{thm:1} gives a result concerning a strip-shaped
(i.e., fixed width; see definition below) subset of $\F$ with gradient $1/4$, by contrast with those of Lemmas~\ref{Fsubset3}
and~\ref{l_strip}, in which the strips have gradient 0 and 1 respectively. Theorem~\ref{thm:2}
gives a construction for a more general subset of $\F$, one which can
lead to a scaled factorial $C(n)$.

\begin{definition}[Strip]
Let $0\leq \gamma\leq 1$, $n_0\in\mathbb N$ and $b > a$. An infinite
sequence $s$ is said to lie in a strip with gradient $\gamma$ if
$$\lfloor\gamma n+a\rfloor\leq s(n) \leq\lfloor\gamma n+b\rfloor$$
for all $n\geq n_0$.
\end{definition}

\subsection{A strip-shaped subset of $\F$}

\begin{theorem} \label{thm:1}
Fix $a>1/4$ and $b>0$, and set
\begin{equation} \label{cd}
c:= \frac{1+4b}{3} , \qquad  d:= \frac{12a-3}{4} . 
\end{equation}
Assume that there exists $f_0 \in \F_{n_0}$ such that $q_0 =Q(f_0)\in \Q_{n_0}$ satisfies$\,$\footnote{For \eqref{cond1}
to be satisfied, $\varepsilon$ cannot be too large: more precisely, we require that $\varepsilon < 5/24$.}
\begin{equation} \label{cond1}
\left( \frac{1}{3}+\ve \right) k \le q_0 (k) \le \left( \frac{3}{4}- \ve \right) k , \qquad  \frac{n_0}{4} \le k \le n_0 ,
\end{equation}
for some $\ve \in (0,d/4b]$ and $n_0\in\mathbb{N}$ such that $n_0\ge 4$ and $\ve n_0 \ge \max \{ 4 c , 2d \}$.
If 
\begin{subequations} \label{cond2}
\begin{align}
f(n)=f_0(n) , & \qquad n\le n_0 ,
\label{cond2a} \\
\frac{n}{4} - a \le f(n) \le \frac{n}{4} + b , & \qquad  n > n_0 ,
\label{cond2b}
\end{align}
\end{subequations}
then there exists $q\in \Q$ such that
\begin{subequations} \label{cond3}
\begin{align}
q(n)=q_0(n) , & \qquad n\le n_0 ,
\label{cond3a} \\
\frac{n}{3} + c \le q(n) \le \frac{3}{4}n -d , & \qquad  n \ge n_0 .
\label{cond3b}
\end{align}
\end{subequations}
\end{theorem}

\begin{proof}
The proof of the bounds \eqref{cond3b} is by induction on $n\ge n_0$.
First, note that, for $n=n_0$, since $\ve n_0 > \max\{c,d\}$, \eqref{cond3b} follows from \eqref{cond1}.
For $n>n_0$, we use the recursion equation
\begin{equation} \label{recursion}
q(n) = q(n-q(n-1)) + f(n)
\end{equation}
to obtain upper and lower bounds on $q(n)$.

\vspace{.1cm}
\noindent\emph{Upper bound.}
If $n-q(n-1) \ge n_0$, then, by \eqref{cond3b} and the inductive hypothesis,
we find 
\[
\begin{aligned}
q(n) &\le \frac{3}{4} \left( n-q(n-1) \right) - d + f(n) \le \frac{3}{4} \left( n - \frac{n-1}{3} - c \right) - d + \frac{n}{4} + b \\
& = \frac{3}{4} n + \frac{1}{4} - \frac{3}{4} c - d + b = \frac{3}{4} n - d + \frac{1 - 3 c + 4 b}{4} = \frac{3}{4}n -d ,
\end{aligned}
\]
because of the definition of $c$ in \eqref{cd}. If $n-q(n-1) < n_0$, then, by using$\,$\footnote{Note that $n-q(n-1) \ge d + n/4$ by the
inductive hypothesis, since $n-1 \ge n_0$, so that $n-q(n-1) \ge n_0/4$.}
\eqref{cond1} for $q(n-q(n-1))$ and the inductive hypothesis for $q(n-1)$, we find
\[
\begin{aligned}
q(n) & \le \left( \frac{3}{4} - \ve \right) \left( n-q(n-1) \right) + f(n) \le \left( \frac{3}{4} - \ve \right) \left( n - \frac{n-1}{3} - c \right) + \frac{n}{4} + b \\
& = \frac{3}{4} n + \frac{1}{4} - \frac{3}{4} c - \frac{2}{3} \ve n - \frac{\ve}{3} + \ve c + b \le \frac{3}{4} n + \frac{1}{4} - \frac{3}{4}c - \frac{4}{3} d  - \frac{\ve}{3} + \ve c + b \\
& =  \frac{3}{4} n - d + \frac{1 - 3c + 4b}{4} - \frac{d - (3c -1) \ve}{3}  = \frac{3}{4}n - d - \frac{d-4b\ve}{3} \le \frac{3}{4}n - d ,
\end{aligned}
\]
where we have used $\ve n > \ve n_0 \ge 2d$ and $4b\ve \le d$, and again, the definition of $c$ in~\eqref{cd}.

\vspace{.1cm}
\noindent\emph{Lower bound.}
We reason in a similar way. If $n-q(n-1) \ge n_0$, we have 
\[
\begin{aligned}
q(n) & \ge \frac{1}{3} \left( n-q(n-1) \right) + c + f(n) \ge \frac{1}{3} \left( n - \frac{3}{4}(n-1) + d \right) + c + \frac{n}{4} - a \\
& = \frac{n}{3} + \frac{1}{4} + \frac{d}{3} + c - a = \frac{n}{3} + c + \frac{3 + 4d- 12a }{12} = \frac{n}{3} + c ,
\end{aligned}
\]
while if $n-q(n-1)<n_0$ we find 
\[
\begin{aligned}
q(n) & \ge \left( \frac{1}{3} + \ve \right) \left( n-q(n-1) \right) + f(n) \ge \left( \frac{1}{3} + \ve \right) \left( n - \frac{3}{4}(n-1) + d \right) + \frac{n}{4} - a  \\
& = \frac{n}{3} + \frac{1}{4} + \frac{d}{3} + \frac{1}{4} \ve n + \frac{3}{4}\ve + \ve d - a \ge \frac{n}{3} + \frac{1}{4} + \frac{d}{3} + c + \frac{3}{4}\ve + \ve d - a \\
& =  \frac{n}{3} + c + \frac{4d + 3 -12a}{12} + \left( \frac{3}{4}\ve + \ve d  \right) \ge \frac{n}{3} + c ,
\end{aligned}
\]
where we have used the definition of $d$ in \eqref{cd} and the condition that $\ve n > \ve n_0 \ge 4c$.
\end{proof}

\begin{remark}
\emph{
As an example take $f(n) = \lfloor (n+2)/4 \rfloor$, which yields $q(n)=\lfloor (n+2)/2 \rfloor$, for $n\le n_0$,
for some $n_0 \in\mathbb{N}$,
so that the bounds \eqref{cond1} are satisfied, with an explicit expression for $\ve$. Theorem \ref{thm:1} shows that if,
for $n>n_0$, we allow $f$ to make moderate jumps, the corresponding sequence $q$ is still well defined. On the other hand,
if the jumps are too large, this is no longer true. For instance, if we take $n_0=4m-1$, with $m\in\mathbb{N}$,
and set $f(n_0+1)=\lfloor 3(n_0+1)/4\rfloor$, then we find $q(n_0+1)=n_0+2$. 
}
\end{remark}

\subsection{A result concerning cones}

As a generalisation of linear cones, we define a cone as a set of sequences $s(n)$
such that $l(n)\leq s(n)\leq u(n)$, where $u(n)\geq l(n)$ for all $n$. 

\begin{theorem} \label{thm:2}
Assume that there exist a real sequence $(a(n))_{n\in\mathbb{N}}$ and an integer $n_0\ge 4$ such that,
for some constants $A,B,C_0\in\mathbb{R}_+$, such that
\begin{equation} \label{AB}
C_0 < \sqrt{\frac{2}{3}} , \qquad \frac{7}{6} + 3 A C_0 \le B \le \frac{2}{C_0} \left( A - \frac{10}{27} \right) ,
\end{equation}
the following holds: 
\begin{enumerate}
\itemsep.2em
\item[(1)] one has
%
\begin{equation} \label{n0n04}
\frac{a(n_0/4)}{n_0/4} < \frac{5}{12(4A+3B)} , \qquad
a(n_0) \ge \frac{1}{2} , \qquad
\frac{1}{3n_0} + \frac{4Aa(n_0)}{n_0} \le \frac{1}{6} \, .
\end{equation}
\item[(2)] the sequence $(a(n))_{n\in\mathbb{N}}$ is increasing for $n\ge n_0/4$ and diverges,
\item[(3)] the sequence $(a(n)/n)_{n\in\mathbb{N}}$ is decreasing for $n\ge n_0/4$ and tends to 0 as $n\to\infty$,
\item[(4)] one has
\vspace{-.3cm}
\begin{equation} \label{anan1-new}
a(n) - a(n') \le C_0 \biggl( 1 - \frac{n'}{n} \biggr) a(n), \qquad n\ge n_0 , \quad \frac{n}{4} \le n'  \le n .
\end{equation}
%
\end{enumerate}
Let $f=(f(n))_{n\in\mathbb{N}}$ be the sequence  such that
\begin{equation} \label{cond4}
\frac{n}{4} - a(n) \le f(n) \le \frac{n}{4} + a(n) , \qquad n \ge n_0 ,
\end{equation}
and assume that $q_{n_0} \in \Q_{n_0}$ exists and satisfies
\begin{equation} \label{cond5}
\frac{n}{3} + 4 A  a(n) \le q(n) \le \frac{3}{4}n - 3 B  a(n) , \qquad \frac{n_0}{4} \le n < n_0 .
\end{equation}
Then, there exists $q=Q(f)\in \Q$ and
\begin{equation} \label{cond6}
\frac{n}{3} + 4 A  a(n) \le q(n) \le \frac{3}{4}n - 3 B  a(n) , \qquad n \ge \frac{n_0}{4} .
\end{equation}
\end{theorem}

\begin{proof}
By induction on $n$. If $n_0/4 \le n < n_0$, the bound \eqref{cond6} follows automatically from \eqref{cond5}.
For $n \ge n_0$, from either \eqref{cond5} or the inductive hypothesis we deduce that
\[
n-q(n-1) \ge n -\frac{3}{4}(n-1) + 3B a(n-1) > \frac{n}{4} \ge \frac{n_0}{4} ,
\]
so that, if $n-q(n-1) < n_0$, then $q(n-q(n-1))$ satisfies the bounds \eqref{cond5} with $n$ replaced with $n-q(n-1)$,
while, if $n-q(n-1) \ge n_0$, then we rely on the inductive hypothesis in order to bound $q(n-q(n-1))$
when looking for bounds on $q(n)$.

To proceed, we use the recursion equation \eqref{recursion} in order to obtain the bounds \eqref{cond6}.
It is convenient to extend the sequence $(a(n))_{n\in \mathbb{N}}$ to a function $a$ defined on $[1,+\infty)$, by setting
\[
a(x) := a(n)+((a(n+1)-a(n))(x-n)
\]
for all $n\in\mathbb{N}$ and all $x\in(n,n+1)$. It is easy to check that the function $a\!:[1,+\infty) \to \mathbb{R}$ satisfies
all the conditions assumed on the sequence $(a(n))_{n\in\,\mathbb{N}}$. In particular, the condition \eqref{anan1-new}
implies that
\begin{equation} \label{anan1}
a(n) - a((1-\delta) n) \le C_0 \delta a(n), \qquad n\ge n_0 , \quad 0 \le \delta \le \delta_0 = \frac{3}{4} .
\end{equation}
Indeed, by writing $(1-\delta)n=m+\xi$, with $m\in\mathbb{N}$ and $\xi\in(0,1)$, we obtain
\[
\begin{aligned}
a(n) - a((1-\delta)n & = a(n) - a(m+\xi) = a(n) - a(m) - (a(m+1)-a(m))\xi \\
& = \left( 1-\xi \right) \left( a(n) - a(m) \right) + \xi \left( a(n) - a(m+1) \right) \\
& \le \left(1-\xi \right) C_0 \biggl( 1 - \frac{m}{n} \biggr) a(n) + \xi C_0 \biggl( 1 - \frac{m+1}{n} \biggr) a(n) \\
& = C_0 \biggl( 1 - \frac{m}{n} \biggr) a(n) - \xi C_0 \frac{a(n)}{n} =
C_0 \biggl( 1 - \frac{m+\xi}{n} \biggr) a(n) \le C_0 \delta a(n) .
\end{aligned}
\]

\vspace{.2cm}
\noindent\emph{Upper bound.}
Fix $n \ge n_0$. From \eqref{cond5}, if $n=n_0$, and from the inductive hypothesis \eqref{cond6}, if $n>n_0$, we obtain
\[
\begin{aligned}
q(n) & \le \frac{3}{4} \left( n-q(n-1) \right) - 3 B  a(n-q(n-1)) + \frac{n}{4} + a(n) \\
& \le \frac{3}{4} \left( n - \frac{n-1}{3} - 4A  a(n-1) \right)
- 3B  a \left( n - \frac{n-1}{3} - 4A  a(n-1) \right) + \frac{n}{4} + a(n) \\
& = \frac{3}{4} n + \frac{1}{4} - 3A  a(n-1) -  3 B  a \left( \frac{2}{3} n + \frac{1}{3} - 4 A  a(n-1) \right) + a(n) ,
\end{aligned}
\]
so that
\[
q(n) \le \frac{3}{4}n - 3B a(n)
\]
if
\begin{equation} \label{need1}
\frac{1}{4} + 3 B  a(n) - 3B  a \left( \frac{2}{3} n + \frac{1}{3} - 4A   a(n-1) \right) \le 3 A  a(n-1) - a(n) .
\end{equation}
Since, by \eqref{anan1},
\begin{equation} \label{anan1bis}
a(n-1)=a(n-1)-a(n) + a(n) \ge a(n) - C_0 \frac{a(n)}{n} ,
\end{equation}
and
\[
a(n) - a \left( \frac{2}{3} n - \frac{1}{3} - 4 A  a(n-1) \right) \le
C_0 \left( \frac{1}{3} + \frac{1}{3n} + \frac{4A a(n-1)}{n} \right) 
a(n) \le \frac{C_0}{2} a(n) ,
\]
where we have used, by condition (3) and the third inequality in \eqref{n0n04},
the fact that
\[
\frac{1}{3n} + \frac{4A a(n-1)}{n} \le \frac{1}{3n} + \frac{4A a(n)}{n} \le
\frac{1}{3n_0} + \frac{4A a(n_0)}{n_0} \le \frac{1}{6} ,
\]
for \eqref{need1} to be satisfied it is enough to require
\[
\frac{1}{4} + \frac{3}{2} B C_0   a(n) \le \left( 3A - 1 - \frac{3AC_0}{n} \right) a(n) ,
\]
that is 
\[
\left( 3A- \frac{3}{2} BC_0  - 1 -\frac{3AC_0 }{n} \right) a(n) \ge \frac{1}{4} .
\]
Since the left hand side in the last equation is increasing in $n$, this means that
\[
\left( 3A- \frac{3}{2} BC_0  - 1 - \frac{3AC_0 }{n_0} \right) a(n_0) \ge \frac{1}{4} ,
\]
which, using, by the first condition in \eqref{AB} and the second and third conditions in \eqref{anan1}, that
\[
\frac{3AC_0 }{n_0} \le \frac{3}{4} \left( \frac{1}{6} - \frac{1}{3n_0} \right) \frac{C_0}{a(n_0)} \le
\frac {C_0}{16 a(n_0)} < \frac{1}{9} , 
\]
is found to be satisfied if we take
\begin{equation} \label{AB1}
3A- \frac{3BC_0}{2} \ge \frac{10}{9} .
\end{equation}

\vspace{.1cm}
\noindent\emph{Lower bound.}
A lower bound for $q(n)$, for $n\ge n_0$, is derived in a similar way. We write
\[
\begin{aligned}
q(n) & \ge \frac{1}{3} \left( n-q(n-1) \right) + 4 A  a(n-q(n-1)) + \frac{n}{4} - a(n) \\
& \ge \frac{1}{3} \! \left( n - \frac{3}{4}(n-1) + 3 B  a(n-1) \right) \!
+ 4 A  a \! \left( n - \frac{3}{4}(n-1) + 3B  a(n-1) \right) \! + \frac{n}{4} - a(n) \\
& = \frac{n}{3} + \frac{1}{4} + B  a(n-1)  + 4 A  a \left( \frac{n}{4}  + \frac{3}{4} + 3  B a(n-1) \right) - a(n) ,
\end{aligned}
\]
so that
\[
q(n) \ge \frac{n}{3} + 4  A a(n)
\]
provided that
\begin{equation} \label{need2}
\frac{1}{4} + B  a(n-1) - a(n) \ge 4 A  a(n) - 4 A  a \left( \frac{n}{4} + \frac{3}{4} +3B   a(n-1) \right) .
\end{equation}
Using \eqref{anan1bis} once more and bounding, by \eqref{anan1},
\[
a(n) - a \left( \frac{n}{4} + \frac{3}{4} + 3  B a(n-1) \right) \le C_0 
\left( \frac{3}{4} - \frac{3}{4n} - \frac{3Ba(n-1)}{n} \right) a(n) \le \frac{3}{4} C_0
\]
then \eqref{need2} follows if we require
\[
3 A C_0  a(n) \le \frac{1}{4} + \left( B  - 1- \frac{B C_0}{n} \right) a(n) ,
\]
that is 
\[
\left( B - 3 AC_0   - 1 -\frac{BC_0}{n} \right) a(n) \ge 0 ,
\]
which, using, by \eqref{AB1}, that
\[
\frac{BC_0}{n} \le \frac{BC_0}{n_0} \le \frac{2A}{n_0} = \frac{4A a(n_0)}{n_0} \frac{1}{2a(n_0)} \le \frac{1}{6}, 
\]
 is satisfied if
 \begin{equation} \label{AB2}
 B - 3 A C_0 \ge \frac{7}{6} .
 \end{equation}
For \eqref{AB1} and \eqref{AB2} to hold at the same time we need $C_0^2<2/3$.
Thus, once $C_0\in(0,\sqrt{2/3})$ has been fixed, we need $A$ and $B$ to satisfy the second condition in \eqref{AB}.
\end{proof}

\begin{remark} \label{rmk4}
\emph{
The first inequality in \eqref{n0n04} is needed in order to have
\begin{equation} \nonumber
\frac{n}{3} + 4 A  a(n) < \frac{3}{4}n - 3 B  a(n)
\end{equation}
in \eqref{cond5}. Of course, in order to allow $q(n_0/4)$ to vary in a non-empty range, we need $n_0$ to be
significantly larger than the minimum value for which the inequality is satisfied.
}
\end{remark}

\begin{remark} \label{rmk2}
\emph{
If we aim to have $A$ and $B$ as small as possible,
we may choose the second condition in \eqref{AB} to hold as a double equality, that is
\begin{equation} \label{AB4}
A = \frac{1}{54} \left( \frac{40+63C_0}{2 - 3 C_0^2} \right) , \qquad
B = \frac{1}{9} \left( \frac{21+20C_0}{2 - 3C_0^2} \right) .
\end{equation}
For instance, 
if $C_0=2/3$ we find $A=41/18 \approx 2.278$ and $B=103/18 \approx 5.722$;
if $C_0=1/2$ we find $A=143/135 \approx 1.059$ and $B=124/25 \approx 2.756$;
if $C_0=1/3$ we find $A=61/90 \approx 0.6778$ and $B=83/45 \approx 1.844$.
For other examples, see Remark \ref{rmk6} below. 
}
\end{remark}

\begin{remark} \label{rmk5}
\emph{
Sequences $(a(n))_{n\in\mathbb{N}}$ which satisfy the hypotheses of Theorem \ref{thm:2} are, for instance, the sequences
which are defined asymptotically as $a(n)=\mu \log n$, with $\mu\in \mathbb{R}_+$,
or $a(n)=\mu n^{\alpha}$, with $\alpha\in(0,1)$ and $\mu\in \mathbb{R}_+$.
In the first case, the value of $C_0$ is arbitrary, while in the latter case we can take$\,$\footnote{We use the
fact that $(1-(1-\delta)^{\alpha})/\delta$ is increasing for $\delta\in(0,1]$ and $\alpha\in(0,1)$.}  
\begin{equation} \label{C02}
C_0 = \frac{4}{3} \left( 1 - \left( \frac{1}{4} \right)^{\alpha} \right) .
\end{equation}
%
The first condition in \eqref{AB}, with $C_0$ of the form \eqref{C02},
implies that, if $a(n)=\mu n^{\alpha}$ for $n$ large enough, then we must require
\[
\alpha < \alpha_0:=-\frac{1}{\log 4} \log \left( 1 - \frac{1}{2} \sqrt{\frac{3}{2}} \right) \approx \frac{2}{3} ,
\]
so that $a(n)$ is $\alpha$-H\"older with exponent $\alpha<\alpha_0$.
Larger values of both $\mu$ and $\alpha$ imply a wider range of variability for both $f$ and $q$;
however \eqref{C02} imposes an upper bound on the value of $\alpha$, and taking large $\mu$ entails
larger values of $n_0$ because of \eqref{n0n04} (see also Remark \ref{rmk6} below),
while we would like to have $n_0$ as small as possible.
}
\end{remark}

\begin{remark} \label{rmk6}
\emph{
Let $(a(n))_{n\in\mathbb{N}}$ be such that $a(n)=\mu n^{\alpha}$, as in Remark \ref{rmk5}.
Then, by the first condition in \eqref{n0n04}, we must require
\[
n_0 \ge 4 \left( \frac{12(4A+3B)}{5} \right)^{1/(1-\alpha)} \mu^{1/(1-\alpha)} .
\]
For $\alpha=1/2$, the last inequality, together with \eqref{AB4} and \eqref{C02}, yields the following values:
\begin{enumerate}
\item$C_0=2/3$;
\item $A=41/18\approx 2.278$ and $B=103/18 \approx 5.722$;
\item $n_0 \ge 1.5910 \times 10^4 \mu^2 $.
\end{enumerate}
If, instead, we take $\alpha=1/3$, we obtain
\begin{enumerate}
\item $C_0 \approx 0.4934$;
\item $A \approx 1.0367$ and $B \approx 2.7012$;
\item $n_0 \ge 6.3769 \times 10^2 \mu^{3/2} $.
\end{enumerate}
We note that the value $\alpha=2/3$ is allowed (since $\alpha_0>2/3$), but the corresponding values
of $A$ and $B$ are very large:
\begin{enumerate}
\item $C_0 \approx 0.8042$;
\item $A \approx 28.0818$ and $B \approx 68.9167$;
\item $n_0 \ge 1.7964 \times 10^9 \mu^3 $.
\end{enumerate}
}
\end{remark}

\section{The Bounding Algorithm}
\label{BAsec}

\subsection{Description}
We now describe an algorithm which is based on an idea used in the
proof of Theorem~1 in~\cite{dilh}. It gives us an investigative tool
for probing finite subsets of $\F$ consisting of sequences $f$ of the form
\begin{equation}
\label{fl_fu}
l(n)\leq f(n)\leq u(n) \qquad \mbox{for all $n\in\mathbb N$,
with $l, u\!:\; \mathbb N\to\mathbb Z$ and $l(1) = u(1) = 0$.}
\end{equation}
The lower and upper bounds on $f$ are given by $l$ and $u$ respectively,
and of interest are examples in which all sequences $f$ obeying these bounds are in $\F$
--- especially if $u(n)-l(n)$ increases without limit as $n$ tends to infinity.

Let $A(n)$ and $B(n)$ denote, respectively, the actual, and the deduced set of values that $q(n)$ can 
take given that $f$ obeys the bounds in~\eqref{fl_fu}.  Here, `deduced' means
`deduced using the Bounding Algorithm' which we are about to introduce.
The aim is to obtain sets $B(n)$ such that $A(n)\subseteq B(n)$ for all $n$ --- see Lemma~\ref{AsubB} ---
and which are easier to construct than $A(n)$.
Of course, if, at some index $n'$, $\max B(n') > n'$ or
$\min B(n') < 1$, then we say that $B$ dies at $n'$, and it is possible,
but, in the light of Lemma~\ref{AsubB}, not guaranteed, that a sequence $q = Q(f)$,
with some $f$ that obeys~\eqref{fl_fu},
also dies. On the other hand, it is always true that $A(n)\subseteq \{1:n\}$ provided that $f\in\F_n$.
Hence, if we have sequences $l$ and $u$ obeying~\eqref{fl_fu} and for which,
for $1\leq i\leq n$,
$\min B(i)\geq 1$ and $\max B(i) \leq i$, then \textit{all} sequences $f$ obeying
$l(i)\leq f(i)\leq u(i)$ are in $\F_n$. 

In order to define the Bounding Algorithm, as a preliminary step we introduce what we call
the Original Bounding Algorithm. We proceed recursively starting from $\Bo(1) := \{1:1\}$.
Let us assume that we know $\Bo(i)$, $1\leq i\leq k$. In order to compute
$\Bo(k+1)$, we use the recursion formula $q(k+1) = q(k+1-q(k)) + f(k+1)$ to
find the set of allowed values of $q(k+1)$ given that $f(k+1) \in \{ l(k+1) : u(k+1)\}$.
We can carry out this computation, given that the set $\Bo(k)$ is known, by the inductive hypothesis,
since then we can find the difference set $\{k+1\} - \Bo(k)$.
Furthermore, each element in this difference set is contained in $\{1:k\}$ and so 
$$q(k+1-q(k))\in\bigcup_{j\in \{k+1\} - \Bo(k)} \Bo(j) $$
is well-defined, since each of the sets in the union is already known.
Thus, we may define $\Bo(k+1)$ as the sumset
\begin{equation}
\label{alg1}
\Bo(k+1) := \{l(k+1) : u(k+1)\} \;\; + \bigcup_{j\in \{k+1\}-\Bo(k)} \Bo(j).
\end{equation}
However, when we use the algorithm in practice, it simplifies the
computation if we further require the sets to consist of consecutive integers.
Therefore, we modify the construction above as follows.
Still proceeding recursively, starting from $B(1) :=\{1:1\}$ and
assuming that $B(i)$ is known for all $i\leq k$, we first consider
\begin{equation}
\label{alg2}
B_*(k+1) := \{l(k+1) : u(k+1)\} \;\; + \bigcup_{j\in \{k+1\}-B(k)} B(j) ,
\end{equation}
then we set $L(k+1) := \min B_*(k+1)$ and $U(k+1) := \max B_*(k+1)$,
and, finally, we define
\begin{equation}
\label{alg3}
B(k+1) := \{ L(k+1) \; : \; U(k+1) \} .
\end{equation}
Note that the set $B_*(k)$ is defined in terms of the sets $B(j)$, with $j<k$, and it is introduced only as an intermediate
step before defining $B(k)$. In fact, $B(k)$ differs from $B_*(k)$ because,
even though the sets $B(j)$ contain, by construction, all the integers between $L(j)$ and $U(j)$ for all $j=1,\ldots,k-1$,
the set $B_*(k)$ may still fail to contain all the integers between $L(k)$ and $U(k)$
--- there may be `gaps' in $B_*(k)$. Whereas it is not a given that equation~\eqref{alg1}
necessarily generates sets of successive integers, it will
nonetheless always be true that $\Bo(k)\subseteq\{L(k):U(k)\}=B(k)$
(see Lemma~\ref{AsubB} below). 
Hence, our assumption gains us simplicity, possibly at the expense of precision.
The simplified version of the algorithm can be rewritten by starting from $L(1)=U(1):=1$ and setting recursively
\begin{equation}
\begin{aligned}
\label{simp_alg1_L}
L(k+1) &= \min\left(\{l(k+1)\;:\; u(k+1)\} \;\; + \bigcup_{j\in \{k+1\}-\{L(k):U(k)\}} L(j)\right)\\
&= l(k+1) + \min_{j\in\{k+1-U(k): k+1-L(k)\}} L(j) \phantom{\bigcup^N}
\end{aligned}
\end{equation}
and
\begin{equation}
\label{simp_alg1_U}
U(k+1) = u(k+1) + \max_{j\in\{k+1-U(k): k+1-L(k)\}} U(j),
\end{equation}
so as to define $B(n):=\{L(n) : U(n)\}$ for all $n\in \mathbb N$.
We refer to this simply as the Bounding Algorithm.

Before we show how it can be used to give insights into 
the behaviour of $Q(f)$ for $f$ in certain sets, we need to establish the following result:
\begin{lemma}
\label{AsubB}
$A(j) \subseteq \Bo(j) \subseteq B(j)$ for all $j\in\mathbb N$.
\end{lemma}
\begin{proof}
We use strong induction.
Let us assume that $A(j)\subseteq \Bo(j) \subseteq B(j)$ for $1\leq j \leq k$.
First of all, note that
\begin{equation}\label{noalg}
A(k+1) \subset \{l(k+1) : u(k+1)\} \;\; + \bigcup_{j\in \{k+1\} - A(k)} A(j)  \qquad \hbox{for all } k \in\mathbb N,
\end{equation}
because once $q(1),....,q(k)$ have been fixed, $q(j)$ with $j=k+1-q(k)$ is also fixed, whereas $A(j)$ contains
all possible values of $q(j)$.
By comparing~\eqref{alg1}, \eqref{alg3} and~\eqref{noalg} we immediately see that
$$\bigcup_{j\in \{k+1\} - A(k)} A(j)  \subseteq \bigcup_{j\in\{k+1\} - \Bo(k)} \Bo(j)
\subseteq\bigcup_{j\in\{k+1\} - B(k)} B(j)$$
since in each relation the right-hand side is a union of at least as many
sets, containing at least as many elements, as the left-hand side.
This implies that $A(k+1) \subset \Bo(k+1) \subset B_*(k+1) \subset B(k+1)$,
where the last inclusion is trivial.
\end{proof}

A short worked example may help to illustrate the workings of the Bounding Algorithm. Let
$l = (0, 0, 0)$ and $u = (0, 1, 2)$.
Then there are six sequences $f$ and these, along with $Q(f)$ corresponding to each, are,
respectively:

\begin{minipage}{0.45\textwidth}
\begin{align*}
f =\; & (0,0,0), (0,0,1), (0,0,2)\\
      & (0,1,0), (0,1,1), (0,1,2)
\end{align*}
\end{minipage}
\begin{minipage}{0.45\textwidth}
\begin{align*}
Q(f) =\; & (1,1,1), (1,1,2), (1,1,3)\\
         & (1,2,1), (1,2,2), (1,2,3).
\end{align*}
\end{minipage}

\medskip
Hence, we have $A(1) = \{1\}$, $A(2)= \{1:2\}$ and $A(3)=\{1:3\}$ but
the Bounding Algorithm gives $B(1) = \{1\}$, $B(2)=\{1:2\}$ and $B(3)= \{1:4\}$. The latter result
implies that $B(4)$ does not exist: we therefore cannot rule out
the possibility that there are sequences $q = Q(f)$ that die at $n=4$.
Interestingly, this example appears to be at odds with Lemma~\ref{Fsubset3},
which shows that any of these six sequences, extended with $f(n)\in\{0:2\}$
for $n>3$, are in fact in $\F$. However, Lemma~\ref{AsubB} implies that there is no
contradiction: the Bounding Algorithm can only ever indicate that $A$ \textit{may} die.

\subsection{The Bounding Algorithm in practice}\label{BAp}
We now present three computations using the Bounding Algorithm.

\subsubsection{Example 1: $l(n) = \lfloor{an - b\sqrt{n}}\rfloor$, $u(n) = \lfloor{an+ b\sqrt{n}}\rfloor$}

In this example, we choose $a = 1/4$, $b = 1/15$, these parameters allowing
us to use the Bounding Algorithm in connection with Theorem~\ref{thm:2}, as
we shall see.

\begin{centering}
\begin{figure}[htbp]
\includegraphics[width=5.5in]{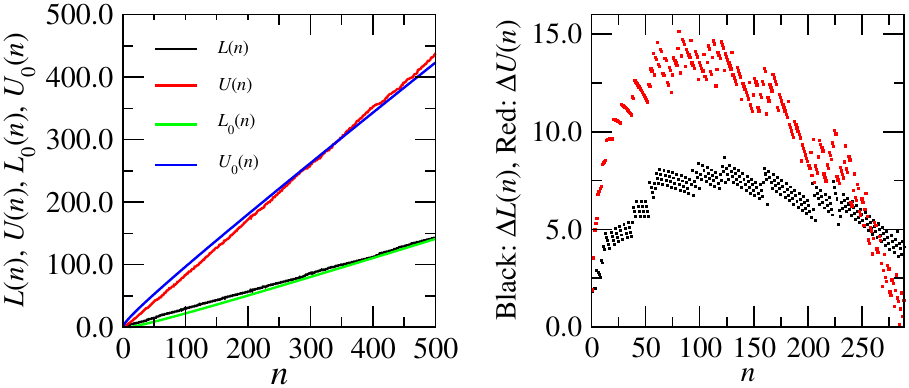}
\caption{Left: plot of $L(n)$, $U(n)$ (black, red respectively), computed using the Bounding
Algorithm, and $L_0(n)$, $U_0(n)$ (green, blue respectively), from~\eqref{boundqbis}.
Right: The differences $\Delta L(n) := L(n)-L_0(n)$ and $\Delta U(n) := U_0(n)-U(n)$,
which should both be positive, whereupon~\eqref{cond5} in Theorem~\ref{thm:2}
is satisfied.}
\label{ex1c}
\end{figure}
\end{centering}

In this example, $l$ and $u$ are not exactly of the form considered in Theorem~\ref{thm:2}, because here we have
$$l(n):=\lfloor an-b\sqrt{n}\rfloor \le f(n) \le u(n):=\lfloor an+b\sqrt{n}\rfloor$$
instead of $$an-b\sqrt{n}\le f(n) \le an+b\sqrt{n}. $$ 
Of course $\lfloor an+b\sqrt{n}\rfloor\le an+b\sqrt{n}$, but, for the lower bound, we have
$an-b\sqrt{n} \ge \lfloor an-b\sqrt{n}\rfloor $. However, we can bound
$$ \lfloor an-b\sqrt{n}\rfloor \ge an- b\sqrt{n} - 1 \ge an - b' \sqrt{n} , \qquad n\ge n_0 $$
provided that
$$ b' \ge b_0 := b + \frac{1}{\sqrt{n_0}}. $$
In that case, if
$$ \lfloor an-b\sqrt{n}\rfloor \le f(n) \le \lfloor an+b\sqrt{n}\rfloor , \qquad n \ge n_0 , $$
then we also have
\begin{equation} \label{boundf}
an-b_0 \sqrt{n}  \le f(n) \le an+ b_0 \sqrt{n} , \qquad n \ge n_0 ,
\end{equation}
with $a = 1/4$ and $b_0=b+(1/\sqrt{n_0})$.

Now, in order to apply Theorem~\ref{thm:2}, we need to check that $q := Q(f)$ satisfies condition~\eqref{cond5}
$$\frac{n}{3} + 4Ab_0 \sqrt{n} \le q(n) \le \frac{3}{4} n - 3B b_0 \sqrt{n} , \qquad \frac{n_0}{4} \le n < n_0,$$
with $C_0 = 2/3$, so that $A = 41/18$, $B =103/18$; and with $n_0$ such that (see Remark \ref{rmk6})
$$n_0 > \left(\frac{24(4A+3B)}{5}\right)^2 \left(b + \frac{1}{\sqrt{n_0}}\right)^2
= 1.591\times 10^4\left(\frac{1}{15} + \frac{1}{\sqrt{n_0}} \right)^2 . $$
This is found to be satisfied for $n_0 > 262.3$, so that, if we choose, say, 
$n_0=17^2=289$, giving $b_0=32/255\approx 0.1255$, we have to check that
$$L_0(n):=\frac{n}{3} + 4Ab_0\, \sqrt{n} \le q(n) \le U_0(n):=\frac{3}{4} n -  3Bb_0 \sqrt{n}, \qquad 72 \le n < 289.$$
Numerically, this amounts to checking that
\begin{equation} \label{boundqbis}
\frac{n}{3} + 1.143\, \sqrt{n} \le q(n) \le \frac{3}{4} n -  2.154 \sqrt{n}, \qquad 72 \le n < 289.
\end{equation}
There are about $\Pi_{n=72}^{289}(U_0(n)-L_0(n)+1)\approx 2.45\times 10^{94}$ 
individual sequences $q$ to check, each of which must obey the
bounds~\eqref{boundqbis}. A direct check of each would clearly be an impossible task.

It turns out, though, that the Bounding Algorithm can be used to show that~\eqref{boundqbis} is 
indeed satisfied --- see
Figure~\ref{ex1c} --- and hence Theorem~\ref{thm:2} applies. We have, in effect,
a computer-assisted proof (since a computer was used to calculate
$L(n)$ and $U(n)$ for $1\leq n\leq 289$) of the following lemma.
\begin{lemma}\label{thm2_example}
For all sequences $f$ obeying
$\lfloor{n/4 - \sqrt{n}/15}\rfloor\leq f(n)\leq \lfloor{n/4 + \sqrt{n}/15}\rfloor$,
$Q(f)\in\Q$ and, for $n\in\mathbb N$, $n/3 + 1.143\, \sqrt{n} \le q(n) \le 3n/4 -  2.154 \sqrt{n}$.
\end{lemma}

\subsubsection{Example 2: $l(n) = \lfloor{a\sqrt{n}}\rfloor$, $u(n) = \lfloor{b\sqrt{n}}\rfloor$}

\begin{centering}
\begin{figure}[htbp]
\includegraphics[width=3.0in]{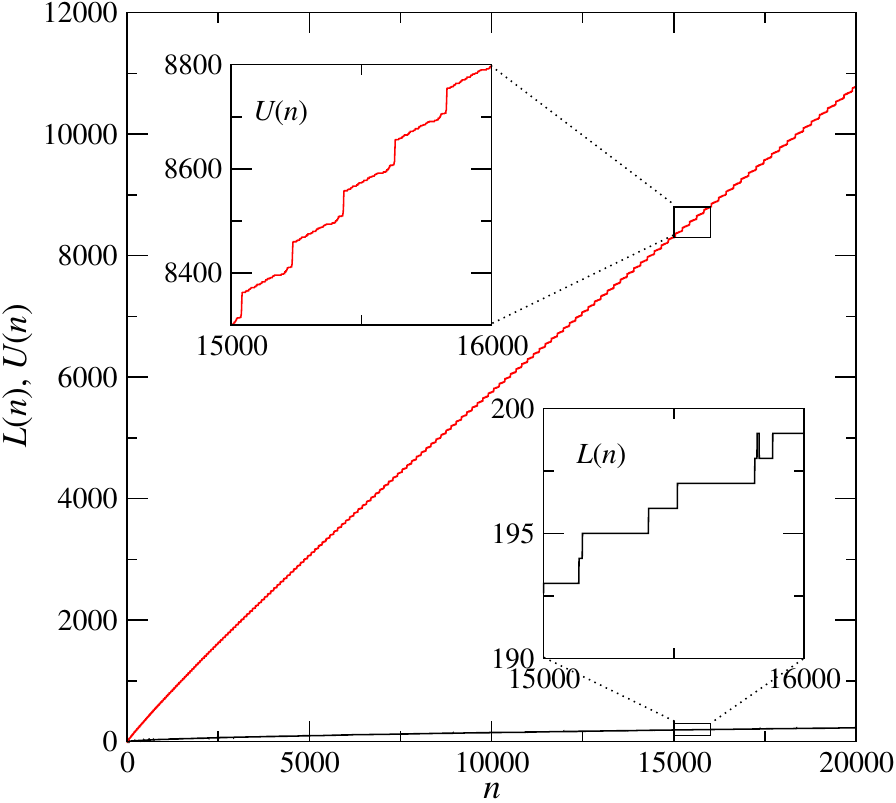}
\caption{Plot for Example 2: $L(n)$ (black) and $U(n)$ (red),
$1\leq n \leq 20000$, computed using the Bounding Algorithm with $l(n) = \lfloor\sqrt{3n/8}\rfloor$,
$u(n) = \lfloor\sqrt{5n/8}\rfloor$. Insets: blow-up of $U(n)$ and $L(n)$ for $15000\leq n \leq 16000$.}
\label{sqrt}
\end{figure}
\end{centering}

As an example not covered by Theorem~\ref{thm:2},
we choose $a = 3/8$, $b = 5/8$ and compute $L(n), U(n)$ for $1\leq n \leq 2^{20}\approx 10^6$.
Over this range, we find that neither $L(n)$ nor $U(n)$ dies, and that 
$L(n)$ and $U(n)$ are well approximated by $L(n)\approx 0.5182 n^{0.5923}$ and
$U(n)\approx 1.4928 n^{0.9176}$ --- see Figure~\ref{sqrt}. Of course, if these approximations hold for all
$n$ then all sequences $l(n)\leq f(n)\leq u(n)$ are in $\F$, although we
currently have no way of knowing if this is true.

\subsubsection{Example 3: $l(n) = \lfloor{an}\rfloor$, $u(n) = \lfloor{bn}\rfloor$}
This is a linear cone. We choose $a = 23/100$ and $b=1/4$ and in this case, $U$ dies at
$n = 845$. A power law fit to the upper bound gives $U(n)\approx 0.5361 n^{1.091}$
which strongly suggests that $U$ dies at some finite $n$.


In summary, the Bounding Algorithm finds application as an investigative tool, since it
can be used to test whether an entire collection of sequences,
rather than just one, is in $\F$. It also leads to proofs, as we saw in
Lemma~\ref{thm2_example}.

\section{Conclusions}

The main objective of this paper was to investigate the set $\F$ of
sequences $f$ such that $Q(f)\in\Q$. To this end, we have proved a selection of results
and a summary of the main ones is given in Table~\ref{summary}.
\begin{table}[h!]
\centering
\begin{tabular}{|l|l|l|l|}\hline
Lemma/Theorem		& $C(n)$	& Shape		& Comments\\ \hline
\ref{monot_not_slow}	& linear	& --- 		& Monotonic but not slow $f$ and $q$\\ \hline
\ref{pwr2}		& exponential	& --- 		& $q'(n)=1,\;\forall n$\\ \hline
\ref{Fsubset3}		& exponential	& Strip		& Horizontal, $C(n)\sim 3^n$\\ \hline
\ref{l_strip}		& exponential	& Strip		& Gradient $=1$, $C(n)\sim \ell^n$, $\ell\geq 3$\\ \hline
\ref{fastest}		& scaled factorial	& ---	& $q'(n)=1,\;\forall n$\\ \hline
\ref{thm:1}		& exponential	& Strip		& Gradient $=1/4$ \\ \hline
\ref{thm:2}		& scaled factorial	& General cone		& Centred on $n/4$ \\ \hline
\ref{thm2_example}	& scaled factorial	& General cone		& Computer-assisted proof \\ \hline
\end{tabular}
\vspace{0.5ex}
\caption{\label{summary} Summary of subsets of $\F$ found and investigated in the paper.}
\end{table}

A variety of approaches was used in the proofs. Most of those in Section~\ref{str_F}
were proved by construction, whereas an analytical approach was adopted in
Section~\ref{further_subsets}. 

Throughout, emphasis has been given to $C(n)$, which
counts the number of sequences $\left(f(i)\right)_{1\leq i\leq n}$ generated
by a particular result. Of particular interest are cases in which
$C(n)\sim \lfloor an\rfloor!$ with $a\in(0,1)$, which we described as displaying `scaled factorial'
growth.

We have also proved some negative results, which we consider to be on an equal
footing with the positive results mentioned above. In particular, we established in
Theorem~\ref{fneg} that any sequence $f$ for which $f(n)<0$ for all $n > n_0\geq 3$
cannot be a member of $\F$. We have also shown (see Lemmas~\ref{cone0}
and~\ref{cone1}) that there are no subsets of $\F$
of the form $0\leq f(n)\leq \lfloor bn\rfloor$ with $b\in (0,1)$, or of the
form $\lfloor(1-\varepsilon)n\rfloor \leq f(n) \leq n$, with
$\varepsilon\in (0, 1)$. It may even be the case that 
linear-cone-shaped subsets of $\F$ cannot exist, but this is just
conjecture. Note, however, that Theorem~\ref{thm:2} does not apply to linear cones,
since these are ruled out by condition (3).

Little has been said in the paper about the original problem underlying this work,
that of the existence, for all $n\in\mathbb N$, of the Hofstadter $q$-sequence~\cite{geb}. This is defined by
$$q_h(n) = q_h(n-q_h(n-1)) + q_h(n-q_h(n-2))\qquad\mbox{with } q_h(1) = q_h(2) = 1.$$
Whereas we can make progress in the case of one nested term, the presence
of two such terms appears to increase the difficulty significantly. 

As a curiosity, we give one simple result implied by Lemma~\ref{Fsubset3}, which is that the sequence
$$q(n) = q(n-q(n-1)) + \left[ q(n-q(n-2))\bmod 3\right] \qquad\mbox{with } q(1)=q(2)=1$$
exists for all $n$.
In fact, computation shows that
$$q(n) = q(n-q(n-1)) +\, \left[ q(n-q(n-2))\bmod M\right] \qquad\mbox{with } q(1)=q(2)=1$$
exists for all $4\leq M\leq 4000$, for $1\leq n\leq2^{22}\approx 4\times 10^6$.
(Clearly, the Hofstadter $q$-sequence is recovered as $M\to\infty$.)

In Section~\ref{BAsec} we introduced the Bounding Algorithm, which, given
that $f$ lies in a cone, that is, $l(n)\leq f(n)\leq u(n)$, allows us
to compute bounds on $q = Q(f)$ of the form $L(n)\leq q(n)\leq U(n)$. We
are unaware of any examples where $L, U$ can be computed from $l, u$
in closed form, but remarkably, the Bounding Algorithm can nonetheless be used, in
conjunction with Theorem~\ref{thm:2}, to give proofs of the existence of
$Q(f)$ with $f$ obeying given bounds. For a concrete example, see Lemma~\ref{thm2_example}.

In~\cite{dilh}, we speculated that equation~\eqref{rec} might turn out to
be an interesting subject for study in its own right,
and this has indeed turned out to be the case.

\end{document}

%% file: header.tex
\documentclass[a4paper,12pt]{amsart}
\usepackage{graphicx}
\usepackage{color}
\usepackage{latexsym}
\usepackage{subfigure}
\usepackage{amscd}
\usepackage{a4wide}
\usepackage{float}
\usepackage{amsmath, amssymb, amsthm, amsfonts, amsgen} 
\usepackage{epstopdf}
\usepackage{mathtools}
\usepackage{enumitem}
\usepackage{etoolbox}
\usepackage{bm}
\usepackage{comment}
\usepackage{pict2e}
\usepackage{relsize}

\graphicspath{ {./} }

\numberwithin{equation}{section}
\newtheorem{definition}{Definition}[section]
\newtheorem{lemma}[definition]{Lemma}
\newtheorem{theorem}[definition]{Theorem}

\newtheorem{corollary}[definition]{Corollary}
\newtheorem{remark}[definition]{Remark}

\newcommand{\ve}{\varepsilon}
\newcommand{\C}{\mathcal{C}}

\newcommand{\F}{\mathcal{F}}

\newcommand{\Q}{\mathcal{Q}}

\newcommand{\Y}{\mathcal{Y}}

\newcommand{\rb}{\rule{0pt}{2.2ex}}
\newcommand{\cstar}{\mathlarger{\star}}

\newcommand{\vp}{\vphantom{3^{3^3}}}
\newcommand{\Bo}{B_{\mathrm{o}}}

%% file: abstract2.tex
\begin{abstract}
We continue work started in~\cite{dilh} concerning
integer sequences $q(n)$, $n\in\mathbb N$, defined by
$q(n) = q(n-q(n-1)) + f(n)$, with $q(1) = 1$. Here,
$f(n)$, with $f(1) = 0$, is a given sequence.
We define $\F$ as the set of
semi-infinite sequences $f$ such that the resulting sequence $q$ exists.
This requires that the term $q(n-q(n-1))$ be defined for all $n$, that is,
$1\leq q(n)\leq n$ applies for all $n\in\mathbb N$.

We use a variety of approaches to probe the structure of $\F$, including
explicit construction, analysis and computer-assisted proof.
\end{abstract}